\documentclass{aptpub}
\usepackage{amsmath}
\usepackage{times}
\usepackage{graphicx}
\usepackage{color}
\usepackage{multirow}
\usepackage{sidecap}

\def\R{{\rm I \mkern-2.5mu \nonscript\mkern-.5mu R}}

\authornames{Ditlevsen and Greenwood} 
\shorttitle{The Morris-Lecar neuron model embeds a leaky
  integrate-and-fire model} 

\begin{document}

\title{The Morris-Lecar neuron model embeds a 
leaky integrate-and-fire model} 

\authorone[University of Copenhagen]{Susanne Ditlevsen} 
\authortwo[University of British Columbia and University of Copenhagen]{Priscilla Greenwood} 

\addressone{Department of Mathematical Sciences, Universitetsparken 5,
DK-2100 Copenhagen {\O}, email: {\tt susanne@math.ku.dk}} 

\addresstwo{Mathematics Annex 1208 2329W East Mall,
Vancouver, BC V6T 1Z4, Canada, email: {\tt pgreenw@math.la.asu.edu}} 

\begin{abstract}
We show that the stochastic Morris-Lecar neuron, in a neighborhood of its
stable point, can be approximated by a two-dimensional
Ornstein-Uhlenbeck (OU)
modulation of a constant circular motion.  The associated radial
OU process is an example of a leaky integrate-and-fire
(LIF) model prior to firing. A new model constructed from a radial OU
process together with a simple firing mechanism based on detailed Morris-Lecar
firing statistics reproduces the Morris-Lecar Interspike Interval (ISI) distribution, and has the
computational advantages of a LIF. The result justifies the large amount of
attention paid to the LIF models. 

\end{abstract}



\section{Introduction}

Much effort has been made to create a realistic but still easily computed
stochastic neuron model, primarily by combining subthreshold dynamics with
firing rules. The result has been a variety of, usually one dimensional, leaky
integrate-and-fire (LIF) descriptions with a fixed membrane potential firing
threshold \cite{Burkitt2006,GerstnerKistler2002,LanskyDitlevsen2008,Lapicque1907}, or with
a rate of firing depending more sensitively on 
membrane potential \cite{Jahn2010,ISI:000237628500003}. These models are useful both for obtaining
analytical results and for ease of simulation.

By contrast, the two-dimensional stochastic Morris-Lecar (ML) neuron model, a
simple cousin to the more detailed Hodgkin-Huxley (HH) model, describes the
dynamics of firing in a way more closely motivated by the biology.
It has been better respected by biologists than the LIF class of models, but
has received little attention owing to the difficulty of mathematical analysis
of this rather complicated stochastic dynamical system.

In Section \ref{sec:StocApprox} of this paper we show that in fact a
LIF model is embedded in the 
ML model as an integral part of it, closely approximating the subthreshold
fluctuations of the ML dynamics.  This result suggests that perhaps
the firing 
pattern of a stochastic ML can be recreated using the embedded LIF together
with a ML stochastic firing mechanism.  We construct such a model in Section
\ref{Sec:firing} and \ref{LIFsection}, and show in Section \ref{ISIs}
that its Interspike Interval (ISI) distribution is similar to that 
of the ML.  Our model, while of the type described in our first paragraph,
combines the realism of the ML with the ease of analysis and
computation of a one dimensional LIF-type model. The work invested in LIF
models is further justified by this new model.

Before we set up our stochastic ML model and write analytical details, let us
have an informal look at how it works. The principal dynamics of the ML,
in the central range of the input current, consist of a stable limit cycle
(Fig. \ref{DeterministicMorrisLecar}A) corresponding to firing, which
encloses a stable fixed point. In
between there loops an unstable limit cycle. The path of the
stochastic model has two quasi-stable patterns
(Fig. \ref{DeterministicMorrisLecar}B). One is 
succesive firings, where the dynamics makes ``large'' noisy circuits
around the stable limit cycle, the other is membrane fluctuations between
spikes, where the dynamics 
makes ``small'' noisy circuits around the fixed point inside the
unstable limit
cycle. The system would continue forever in one of these two patterns
were it not for the noise which causes switching from firing
to subthreshold fluctuations and back again at random times when
the dynamics cross the unstable limit
cycle. Our analysis will show that the dynamics between spikes, of random
cycling inside the unstable limit cycle followed by crossing to the
stable limit cycle outside it, can be identified with the sample path
behavior of a two-dimensional Ornstein-Uhlenbeck (OU) process times a
rotation. 

A main ingredient in our result is the stochastic
dynamical phenomenon that oscillations which damp to a fixed point in
a deterministic system will be sustained by the stochasticity in a
corresponding stochastic system. Damped oscillations in a
two-dimensional system are signalled by a local linear structure
defined by a matrix having a pair of conjugate complex eigenvalues
with negative real part. A corresponding stochastic system will not
damp, being prevented by the noise. Instead, a quasi-stationary
stochastic process is set up, which cycles in a random pattern around
the fixed point. Using recent results of \cite{BaxendaleGreenwood2010}
we are able to identify, approximately, this stochastic process which
is part 
of the subthreshold dynamics of the ML. Up to a fixed linear 
transformation, the approximating process is the product of
a steady fast rotation with a 
two-dimensional OU process. The identification allows us to cement in
place the correspondance, for a particular set of model parameters,
a particular LIF model as the appropriate subthreshold phase between ML
firings.



\section{The Morris-Lecar model}

There exists a large variety of modeling
approaches to the generation of spike trains in
   neurons (see
   e.g. \cite{DayanAbbott2001,GerstnerKistler2002,BookIz}). Most
   famous is 
   the Hodgkin-Huxley (HH) model \cite{HH1952} 
consisting of four coupled 
differential equations, one for the membrane voltage, and three
equations describing the gating variables that model the
voltage-dependent sodium and potassium channels. A large amount of
research effort is currently directed towards 
understanding how neural coding carries information through nervous
systems. Basic to the subject is how single neurons transmit
information. As in any modeling
effort, we must ignore or summarize details and focus on what, we hope,
are a few essential aspects. The ML model \cite{MorrisLecar1981} has
often been used as a good, qualitatively quite accurate,
two-dimensional model of neuronal spiking. It is a conductance-based
model like the HH model, introduced to explain the dynamics of the barnacle muscle
fiber. The original ML model was three-dimensional, including a fast
responding voltage-sensitive Ca$^{2+}$ conductance, and a delayed
voltage-dependent K$^+$ conductance for recovery. To justify the
two-dimensional version, one uses that the Ca$^{2+}$ activation moves on a
much faster time scale than the other variables, and can conveniently
be treated as an instantaneous variable, by replacing it by its
steady-state value given the other variables.  

\begin{table}
\begin{center}
\caption{ \label{MLparameters} Variables and parameter values used in
  the Morris-Lecar model}
\begin{tabular}{c@{ = }r l}
\hline
\multicolumn{2}{l}{$V(t)$  [mV]} & Membrane voltage\\
\multicolumn{2}{l}{$W(t)$  [1]} & Normalized K$^{+}$ conductance \\
\multicolumn{2}{l}{$t$  [ms]} & Time \\
\hline
$V_1$  &  -1.2  mV & Scaling parameter\\
$V_2$  &  18  mV & Scaling parameter\\
$V_3$  &   2  mV & Scaling parameter\\
$V_4$  &  30  mV & Scaling parameter\\
$g_{\text{Ca}}$ &   4.4  $\mu$S/cm$^2$ & Maximal conductance
associated with Ca$^{2+}$ current\\
$g_K$  &   8  $\mu$S/cm$^2$ & Maximal conductance
associated with K$^{+}$ current\\
$g_L$  &   2  $\mu$S/cm$^2$ & Conductance
associated with leak current\\
$V_{\text{Ca}}$ & 120  mV & Reversal potential for Ca$^{2+}$ current\\
$V_K$  & -84  mV & Reversal potential for K$^+$ current\\
$V_L$  & -60  mV & Reversal potential for leak current\\
$C$   &  20  $\mu$F/cm$^2$ & Membrane capacitance\\
$\phi$ &   0.04 1/ms & Rate scaling parameter\\
$I$   &  90  $\mu$A/cm$^2$ & Input current\\
\hline
\end{tabular}
\end{center}
\end{table}

The parameter values in our computations were chosen from
\cite{RinzelErmentrout1998,TatenoPakdaman2004}, and are given in Table
\ref{MLparameters} together with the interpretation of variables and
parameters. The variable $V_t$ represents the membrane potential of
the neuron at time $t$, and $W_t$ represents the normalized
conductance of the K$^+$ current. This is a variable between 0 and 1,
and could be interpreted as the probability that a K$^+$ ion channel
is open at time $t$. The non-linear model equations are 
\begin{eqnarray}
\label{dV}
dV_t &=& \frac{1}{C} \left (-g_{Ca}  m_{\infty}(V_t) (V_t - V_{Ca})
         -g_K W_t (V_t - V_K)-g_L  (V_t - V_L)+I
         \right ) dt, \hspace{5mm} \\
\label{dW}
dW_t &=& \left ( \alpha (V_t)(1-W_t)-\beta (V_t) W_t \right ) dt,
\end{eqnarray}
with the auxiliary functions given by
\begin{eqnarray}
\label{m}
m_{\infty}(v) &=& \frac{1}{2} \left (1 + \tanh \left (\frac{v-V_1}{V_2}
  \right ) \right ), \\
  \alpha (v) &=& \frac{1}{2} \phi \cosh \left (\frac{v-V_3}{2V_4} \right
  ) \left (1 + \tanh  \left (\frac{v-V_3}{V_4} \right ) \right ),\\
\label{beta}
  \beta  (v) &=&  \frac{1}{2}\phi \cosh \left (\frac{v-V_3}{2V_4} \right
  ) \left (1 - \tanh
  \left (\frac{v-V_3}{V_4} \right )\right ).
\end{eqnarray}
Equation {\eqref{dV} describing the dynamics of $V_t$ contains four terms,
  corresponding to Ca$^{2+}$ current, K$^+$ current, a general leak
  current, and the input current $I$. The functions
  $\alpha (\cdot )$ and $\beta (\cdot )$ model the rates of opening
  and closing, respectively, of the K$^+$ ion channels. The function
  $m_{\infty}(\cdot)$ represents the equilibrium value of the
  normalized  Ca$^{2+}$ conductance for a given value of the membrane
  potential. 

In Fig. \ref{DeterministicMorrisLecar}A the phase-state of the
model is plotted. The system has two stable at\-tractors; a stable fixed
point corresponding to quiescence of the neuron, and a stable limit
cycle corresponding to repetitive firing. In between the two
attractors is an unstable limit cycle, which splits the state space into
two parts from either of which the deterministic process cannot escape, once
trapped there. 

\begin{figure}[h]
\begin{picture}(400,190)
\put(35,180){{\large A}}
\put(235,180){{\large B}}
\put(0,0){\includegraphics[height=6.2cm]{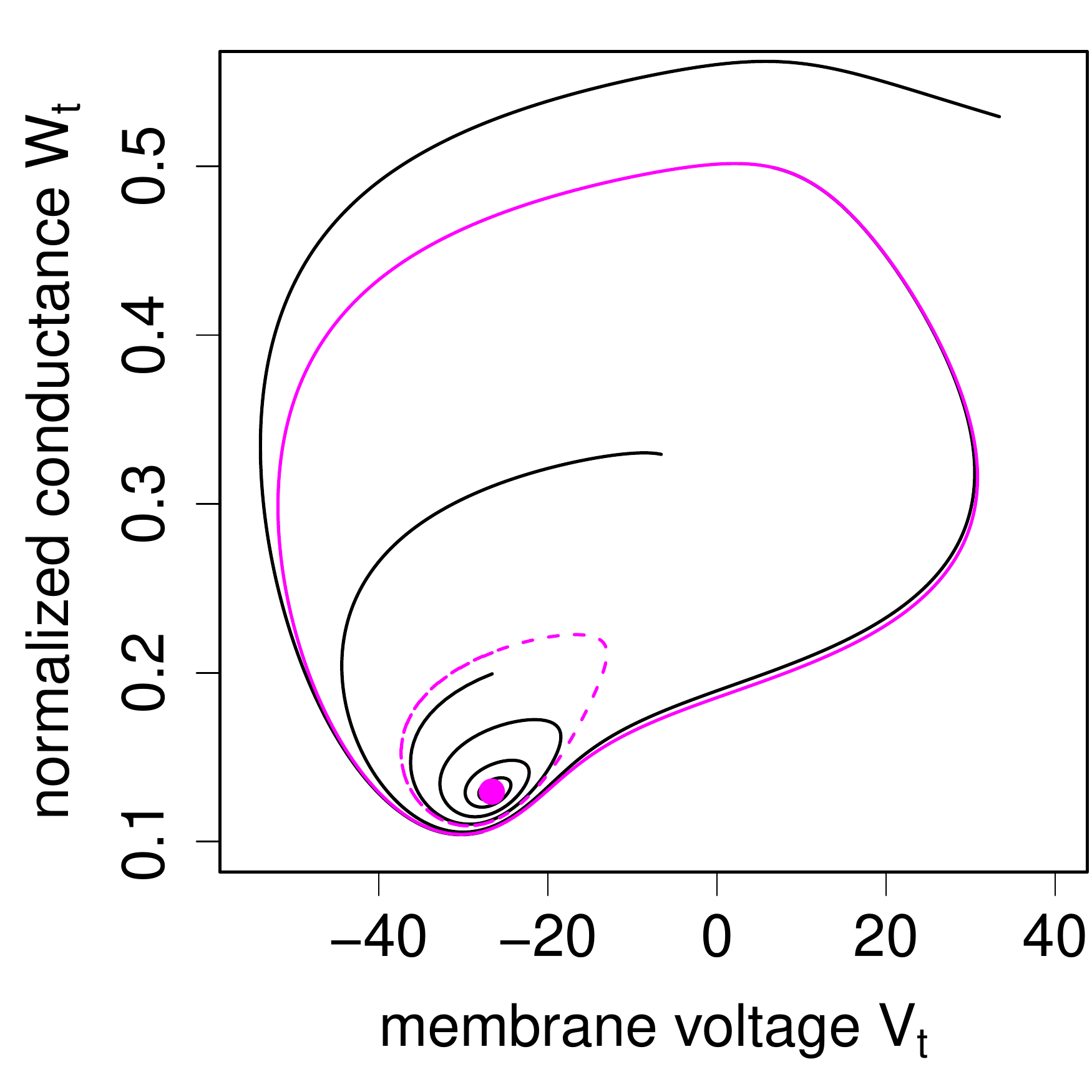}}
\put(200,0){\includegraphics[height=6.2cm]{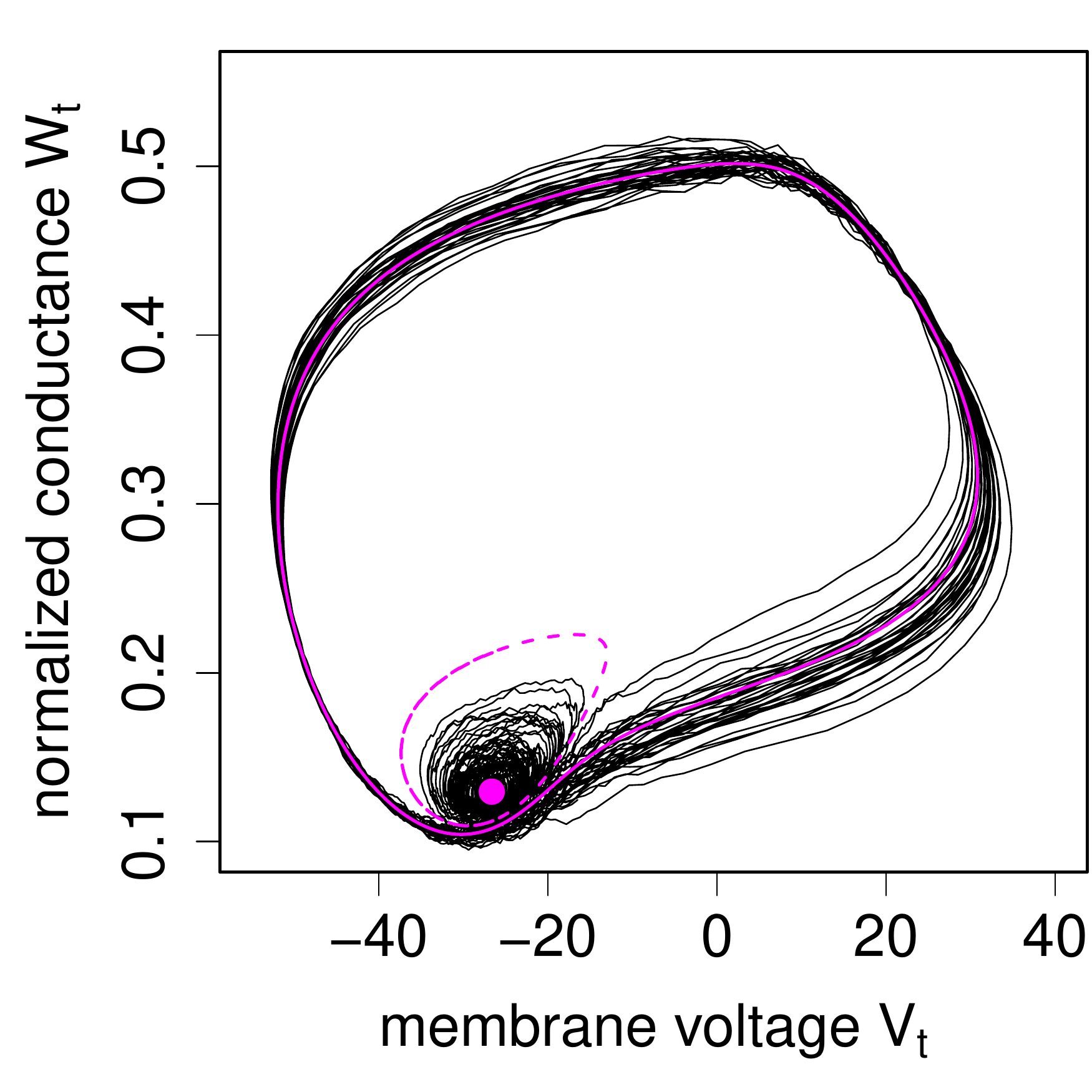}}
\end{picture}
\caption{ \label{DeterministicMorrisLecar} Phase-state plots of the
  normalized conductance $W_t$ against membrane voltage $V_t$. The
  full drawn magenta curve is a
  stable limit cycle, the dashed magenta curve is an unstable limit cycle,
and the magenta point is a stable fixed point. Black curves are sample
trajectories. Panel A: model without noise, \eqref{dV}--\eqref{beta}. If the process is started
between the stable and the unstable limit cycle, or outside the stable
limit cycle, the solution is seen
to spiral out, respectively in, towards the stable limit cycle,
corresponding to repetitive firing of the neuron. If the process is
started inside the unstable limit cycle, the solution spirals into the
stable fixed point, corresponding to subthreshold fluctuations of the neuron. Note
that three trajectories are plotted. Panel B: model with noise,
\eqref{dV}, \eqref{m}--\eqref{beta} and
\eqref{dWsigma}, $\sigma^* =0.05$. Only one trajectory is plotted, and the solution is 
seen to switch between periods of firing and quiescence. }
\end{figure}

\subsection{The stochastic Morris-Lecar model with channel noise}

It has long been known that the opening
and closing of ion channels is an important part of neuron
function. Channel activity is summarized, even in the comparatively
detailed HH model, by potential dependent
averages. However, it has become apparent that the stochastic
nature of ion channels must be explicitly modeled if we are to capture
essential features of neuron dynamics. Changes in
the states of channels cannot be tracked explicitly because of their
vast number. Hence, it is useful to model the role of ion channels as a
stochastic process, $W_t$, the proportion of channels open at time
$t$. We therefore add channel noise by changing the ordinary
differential equation system \eqref{dV} -- \eqref{beta}, to a
stochastic differential equation system,
replacing the conductance equation \eqref{dW} by 
\begin{eqnarray}
\label{dWsigma1}
dW_t &=& \left ( \alpha (V_t)(1-W_t)-\beta (V_t) W_t \right ) dt
         + h (V_t, W_t)  dB_t,
\end{eqnarray}
where $B_t$ is a standard Wiener process, and the function $h (
\cdot )$ has to be chosen. 

The diffusion coefficient $h (\cdot)$ in \eqref{dWsigma1} should
be based on the drift 
coefficient which gives the rate of change of fraction of open ion channels
due to random openings and closings. A natural choice of the function
$h (\cdot)$, following the diffusion approximation of
\cite{Kurtz1978}, would be the  
square root of the sum of the two rates in the drift coefficient, times a
factor $1/\sqrt N$ where $N$ is the number of ion channels involved. However, this
choice has the problem that it is not zero when all the channels are closed, and
the resulting \eqref{dWsigma1} would produce negative solutions with
positive probability. 
To avoid this difficulty, for fixed $V_t$ we let $W_t$ be a Jacobi diffusion. In
fact, in the class of Pearson diffusions \cite{PearsonDiffusions},
i.e. one-dimensional diffusions with linear drift, and with $h^2
(\cdot)$ a polynomial of at most degree two, this is the
only bounded diffusion. Living
on $(0,1)$, it has the form   
\begin{eqnarray}\label{SDEjacobi}
dX_t &=& -\theta \left ( X_t - \mu \right ) dt
         + \gamma \sqrt{2\theta X_t(1-X_t)}  dB_t
\end{eqnarray}
where $\theta > 0$ and $\mu \in (0,1)$. 
It is named for the eigenfunctions of the generator, which are the Jacobi
polynomials. It is ergodic provided that 
$\gamma^2 \leq \min (\mu, (1-\mu))$, 
and its stationary distribution is the Beta distribution with shape
parameters $\mu/\gamma^2$ and $(1-\mu)/\gamma^2$. It has mean $\mu$
and variance $\gamma^2 \mu (1-\mu)/(1+\gamma^2)$. In our case, because the
diffusion coefficient in \eqref{SDEjacobi} should be of the same order
as the one given by the Kurtz approximation
\cite{Kurtz1978}, $\gamma$ is 
proportional to $1/\sqrt N$. 

By equating the drift terms in  
\eqref{dWsigma1} and \eqref{SDEjacobi}, we have $\theta = \alpha (V_t)
+ \beta (V_t)$ and 
$\mu = \alpha (V_t)/$ $ (\alpha (V_t) + \beta (V_t))$. So for fixed $V_t$,
with $h^2 (V_t, W_t) = \gamma^2 2(\alpha (V_t) + \beta (V_t))
W_t(1-W_t)$, where $\gamma^2$ is constrained by $\gamma^2 (\alpha (V_t) +
\beta (V_t))\leq \min (\alpha (V_t), \beta (V_t))$,  
also \eqref{dWsigma1} will stay bounded in $(0,1)$. Since
$\alpha (V_t)$ and $\beta (V_t)$ are strictly positive, we can put $\gamma^2 =
(\sigma^*)^2 \alpha (V_t) \beta (V_t) /(\alpha (V_t) + \beta (V_t))^2$, with
$\sigma^* \in (0,1]$, and specify the conductance equation \eqref{dWsigma1} as
\begin{eqnarray}
\label{dWsigma}
dW_t &=& \left ( \alpha (V_t)(1-W_t)-\beta (V_t) W_t \right ) dt
         + \sigma^* \sqrt{2 \frac{\alpha (V_t) \beta (V_t)}{\alpha
             (V_t) + \beta (V_t)} W_t(1-W_t) } dB_t. 
\end{eqnarray}
In the next Section we compute the equilibrium point
$(V_{\mbox{eq}},W_{\mbox{eq}})$ of the 
system \eqref{dV}--\eqref{beta} for the chosen parameters. By equating the
diffusion coefficient as it would occur in the diffusion approximation
of \cite{Kurtz1978} with the one in \eqref{dWsigma} at
$(V_{\mbox{eq}},W_{\mbox{eq}})$ we will obtain $\sigma^*$ in terms of
$1/\sqrt{N}$, where $N$ is the number of channels involved.

It can be shown by a coupling argument that also for varying $V_t$
will $W_t$ given by \eqref{dWsigma} stay
bounded in $(0,1)$, since $V_t$ is bounded once it
is started inside some interval \cite{DitlevsenJacobsen2011}.


In Fig. \ref{StochasticMorrisLecar}, the model defined by \eqref{dV},
\eqref{m}--\eqref{beta} and \eqref{dWsigma} is 
simulated for different values of $\sigma^*$, where these can be
thought of as corresponding to different total numbers of ion 
channels. 

\begin{figure}[h!]
\centerline{\includegraphics[height=6cm]{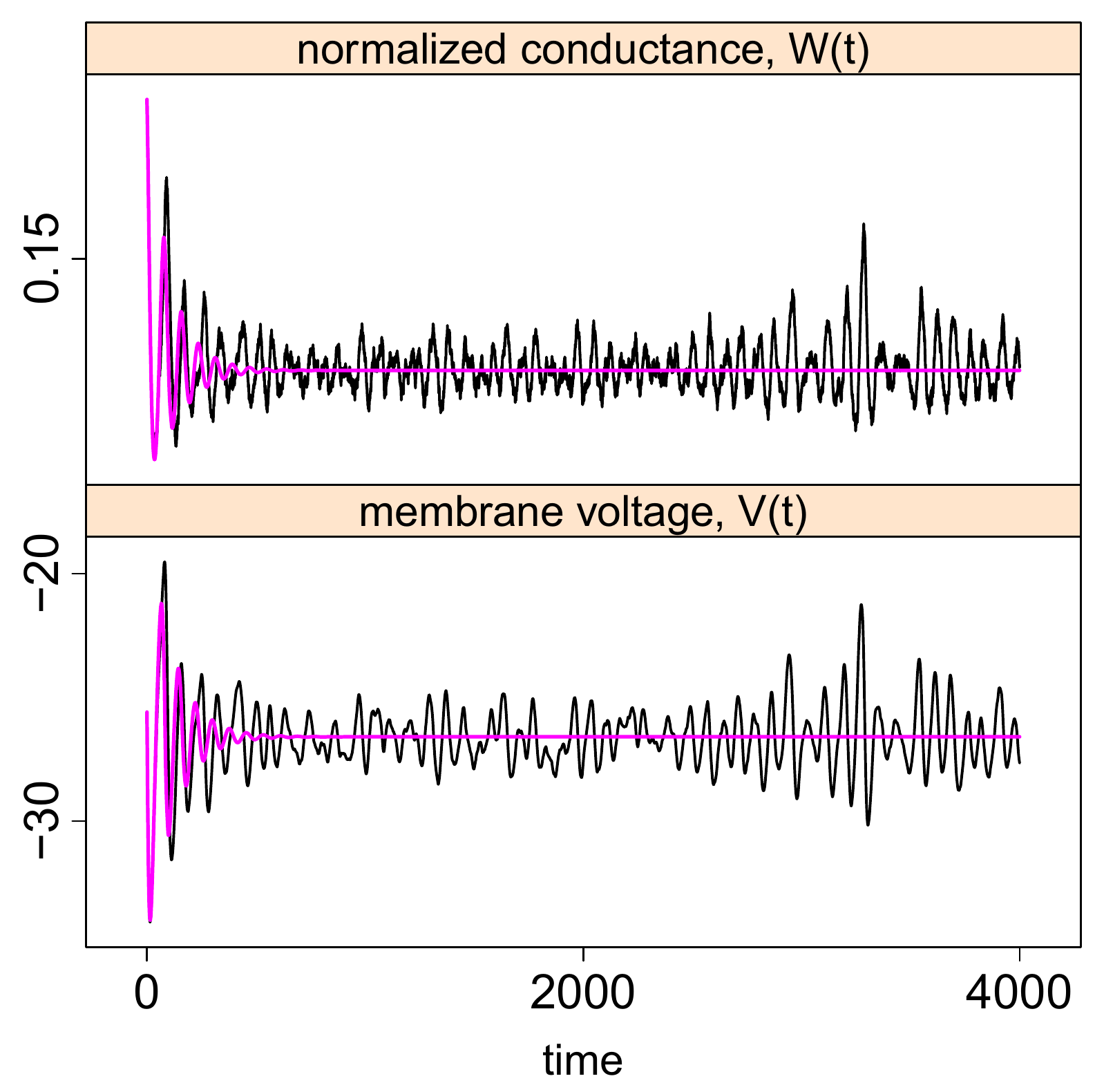}
\includegraphics[height=6cm]{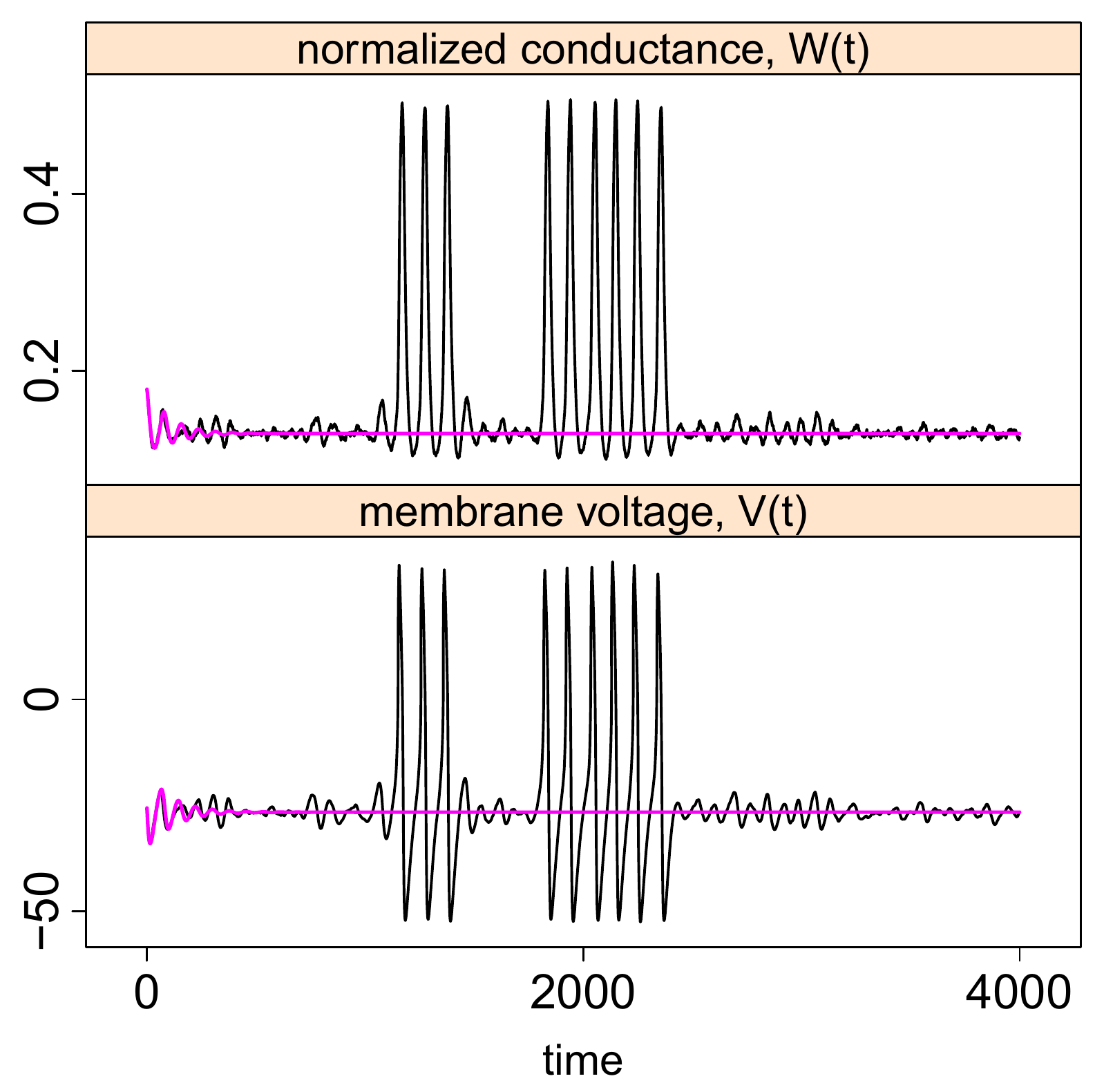}}
\centerline{\includegraphics[height=6cm]{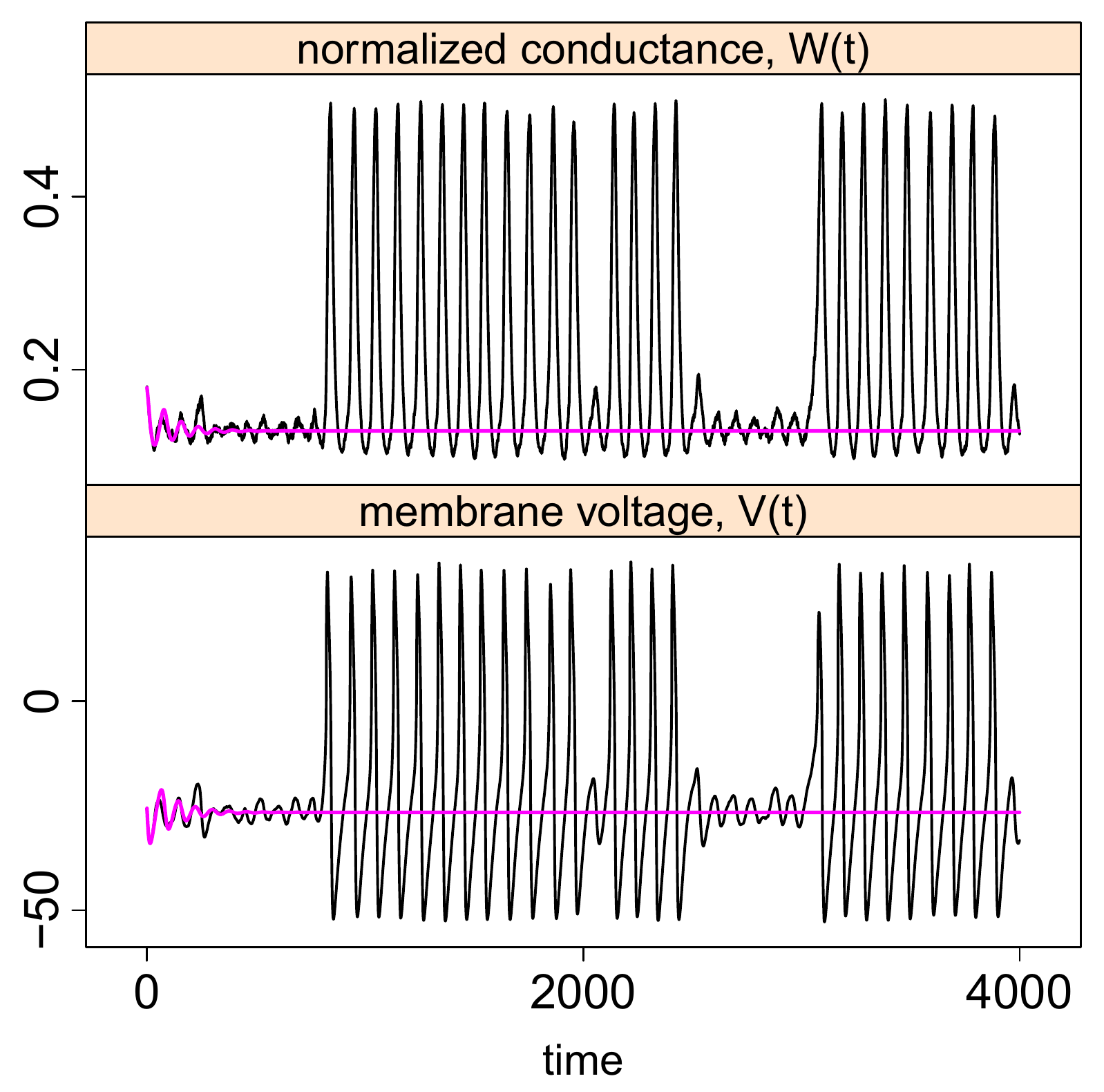}
\includegraphics[height=6cm]{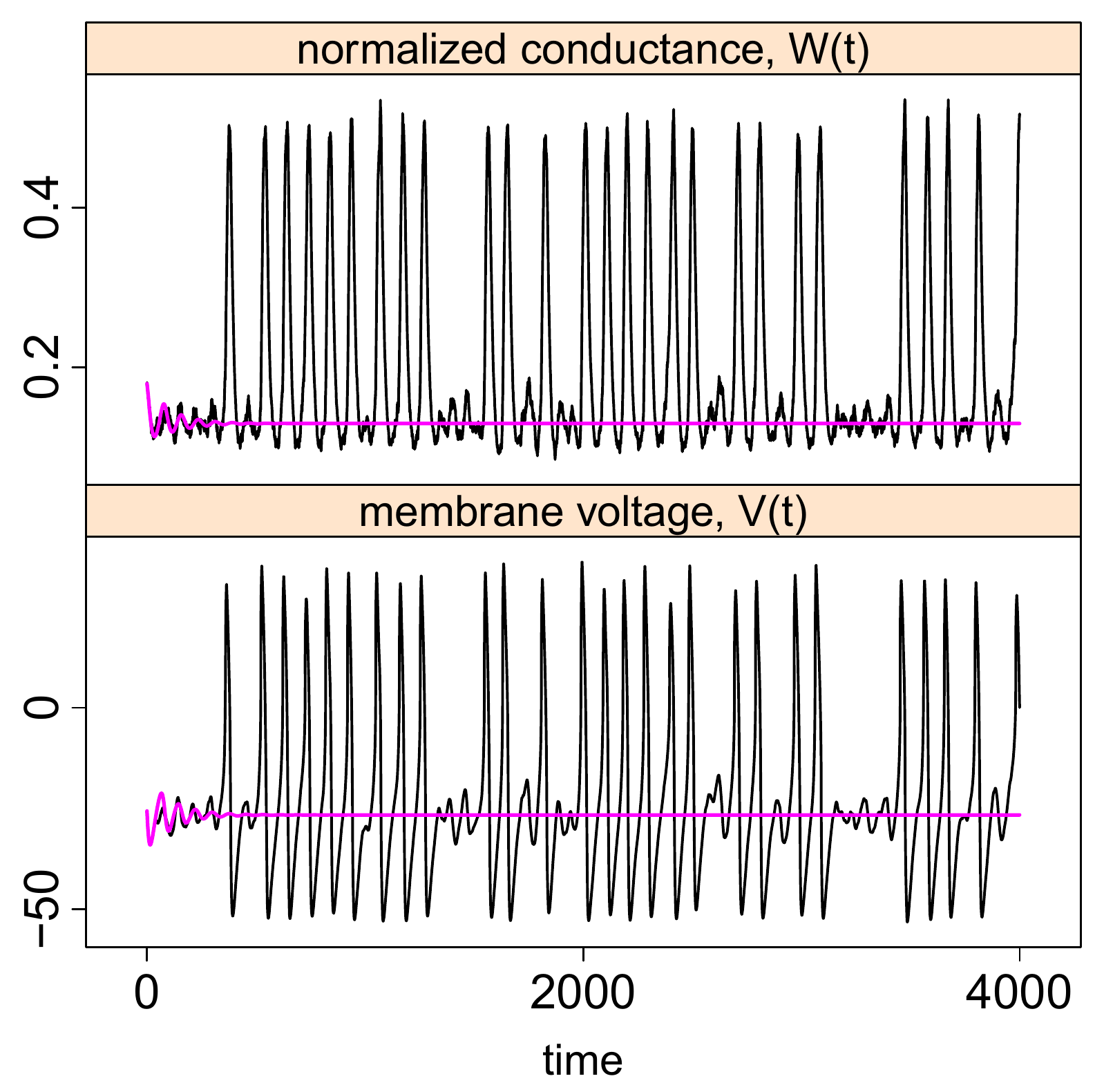}}
\caption{ \label{StochasticMorrisLecar} Time series plots (black
  curves)  of the
  stochastic Morris-Lecar model for different noise levels started
  inside the unstable limit cycle, but not at the fixed point. Upper left:
  $\sigma^* =0.02$, upper right: $\sigma^* =0.03$, lower left: $\sigma^*
  =0.05$, lower right: $\sigma^*
  =0.1$. Note different scales, in the upper left panel there is no
  firing. The magenta curves are the deterministic model, $\sigma^*=0$. }
\end{figure}

\section{The linear approximation of the stochastic Morris-Lecar
  during quiescence} 

To identify the process of subthreshold oscillations, i.e. the dynamics
close to the stable fixed point between firings, we analyze the
linearized system around this point. Consider the system
\begin{eqnarray*}
dV_t &=& f(V_t, W_t)dt, \\
dW_t &=& g(V_t, W_t)dt + h(V_t, W_t)dB_t,
\end{eqnarray*}
where the functions $f(\cdot)$, $g(\cdot)$ and $h(\cdot)$ are given by
\eqref{dV}, \eqref{m}--\eqref{beta} and \eqref{dWsigma}.

For the chosen parameter values given in Table 
\ref{MLparameters}, the deterministic system, obtained for
$h(\cdot)=0$, has a unique locally stable equilibrium point 
$(V_{\text{eq}},W_{\text{eq}})$ given by
\begin{eqnarray*}
W_{\text{eq}}(V_{\text{eq}}) &=& \frac{\alpha (V_{\text{eq}})}{\alpha (V_{\text{eq}})
+ \beta (V_{\text{eq}})} \, = \, \frac{1}{2}
\left (1 + \tanh  \left (\frac{V_{\text{eq}}-V_3}{V_4} \right ) \right )
\end{eqnarray*}
and $V_{\text{eq}}$ is the solution to the equation
$f(V_{\text{eq}},W_{\text{eq}}(V_{\text{eq}}))=0$, which cannot be
solved analytically, but can be
found numerically. The input current value $I=90\mu$A/cm$^2$ is a
typical value well 
inside the range of $I$ where the deterministic dynamics has a
stable limit point inside an unstable limit cycle as shown in Fig.
\ref{DeterministicMorrisLecar}A. The equilibrium point for
$I=90\mu$A/cm$^2$ is
\begin{eqnarray*}
(V_{\text{eq}},W_{\text{eq}}) &=& (-26.6 \text{ mV }, \,  0.129).
\end{eqnarray*}
In terms of the centered variables
\begin{eqnarray*}
X_t^{(1)} \, = \, V_t - V_{\text{eq}} &,&
X_t^{(2)} \, = \, W_t - W_{\text{eq}}
\end{eqnarray*}
the system becomes
\begin{eqnarray}
\label{dXt}
dX_t^{(1)} &=& f \left (
  (X_t^{(1)}\! +V_{\text{eq}}),(X_t^{(2)}\! +W_{\text{eq}}) \right )dt + 0 \cdot 
dB_t^{(1)} \nonumber \\
& = & f^* (X_t^{(1)},X_t^{(2)})dt, \\
\label{dYt}
dX_t^{(2)} &=& g \left ( (X_t^{(1)}\! + V_{\text{eq}}),(X_t^{(2)} \!
  +W_{\text{eq}}) \right )dt +  
h \left ( (X_t^{(1)}\! +V_{\text{eq}}),(X_t^{(2)}\! +W_{\text{eq}})
\right )dB_t^{(2)} \hspace{3mm} \nonumber \\
& = & g^* (X_t^{(1)},X_t^{(2)}) dt + h^* (X_t^{(1)},X_t^{(2)}) dB_t^{(2)}.
\end{eqnarray}
We write $X_t = (X_t^{(1)} \quad X_t^{(2)})^T$ and $B_t
= (B_t^{(1)} \quad B_t^{(2)})^T$, where
$^T$ denotes transposition. Note that $B_t^{(1)}$ does not enter the
dynamics, but is introduced to ease the matrix notation, as will be clear in the
following. When the noise is small and the process $X_t$ is started
near the equilibrium point,  $x=(0,0)$, we expect the dynamics to
concentrate around 
the equilibrium point. A local approximation is obtained by
linearizing \eqref{dXt}--\eqref{dYt} around $(0,0)$. The
diffusion term is approximated by setting $X_t^{(1)} = X_t^{(2)} =0$ in the diffusion
coefficients. The linearized system is 
\begin{eqnarray}
\label{linearsystem}
d X_t &=&
\mathbf{M} X_t dt +
\mathbf{G}  dB_t,
\end{eqnarray}
where
\begin{eqnarray*}
\mathbf{M} &=& \begin{pmatrix} m_{11} & m_{12} \\ m_{21} & m_{22} \end{pmatrix}
 \, = \, \left . \begin{pmatrix} \frac{\partial f^*}{\partial x_1} &
  \frac{\partial f^*}{\partial x_2}\\ \frac{\partial g^*}{\partial x_1} &
  \frac{\partial g^*}{\partial x_2} \end{pmatrix} \right
|_{(x_1,x_2)=(0,0)} \! = \,  \begin{pmatrix}  0.0258 &  -22.961 \\
0.000335 & -0.0446 \end{pmatrix},
\end{eqnarray*}
using the parameter values in Table \ref{MLparameters}, and
\begin{eqnarray}
\label{G}
\mathbf{G} &=& \begin{pmatrix} 0&0 \\[2mm] 0&\sigma^*
  \sqrt{2(\alpha(V_{\text{eq}})+\beta(V_{\text{eq}}))}
  (1-W_{\text{eq}})W_{\text{eq}}\end{pmatrix} 
\, =
\, \begin{pmatrix} 0&0 \\[2mm] 0& \sigma \end{pmatrix},
\end{eqnarray}
where $\sigma = 0.034 \sigma^*$. 
By evaluating the diffusion approximation of \cite{Kurtz1978} at
$(V_{\text{eq}},W_{\text{eq}})$ and
equating to the above we obtain $\sigma^* =
1/\sqrt{W_{\text{eq}}(1-W_{\text{eq}})N} \approx 3/\sqrt{N}$. 
In the Appendix the matrix $\mathbf{M}$ is detailed.
Solutions of \eqref{linearsystem} with $\mathbf{G} =0$ are given
in terms of the eigenvalues of $\mathbf{M}$ which are complex
conjugates and given by
\begin{eqnarray*}
-\lambda \pm \omega i &=& -0.0094 \pm 0.0803 i
\end{eqnarray*}
where $\lambda = -\text{tr} (\mathbf{M})/2$, $\omega^2 = |
  \lambda^2 - \det (\mathbf{M})|$ and $i=\sqrt{-1}$.
Thus, near the equilibrium point the solution of \eqref{linearsystem},
with $\sigma =0$, is
\begin{eqnarray}
\label{approxsolution}
X_t &=&
\mathbf{C} \begin{pmatrix} \cos \omega t\\ \sin \omega
  t   \end{pmatrix} e^{-\lambda t},
\end{eqnarray}
where $\mathbf{C}$ contains the initial conditions
\begin{eqnarray*}
\mathbf{C} &=& \begin{pmatrix} x_0 &  (m_{12} y_0 + (m_{11}
    + \lambda) x_0)/\omega \\ y_0 & (m_{21} x_0 - (m_{11}
    + \lambda) y_0)/\omega \end{pmatrix}.
\end{eqnarray*}
In Fig. \ref{DeterministicApprox} the solution of the deterministic model,
\eqref{dV}--\eqref{beta} with $\sigma =0$, is compared to the linear approximation
\eqref{approxsolution}.

\begin{SCfigure}
  \centering
\includegraphics[width=8cm]{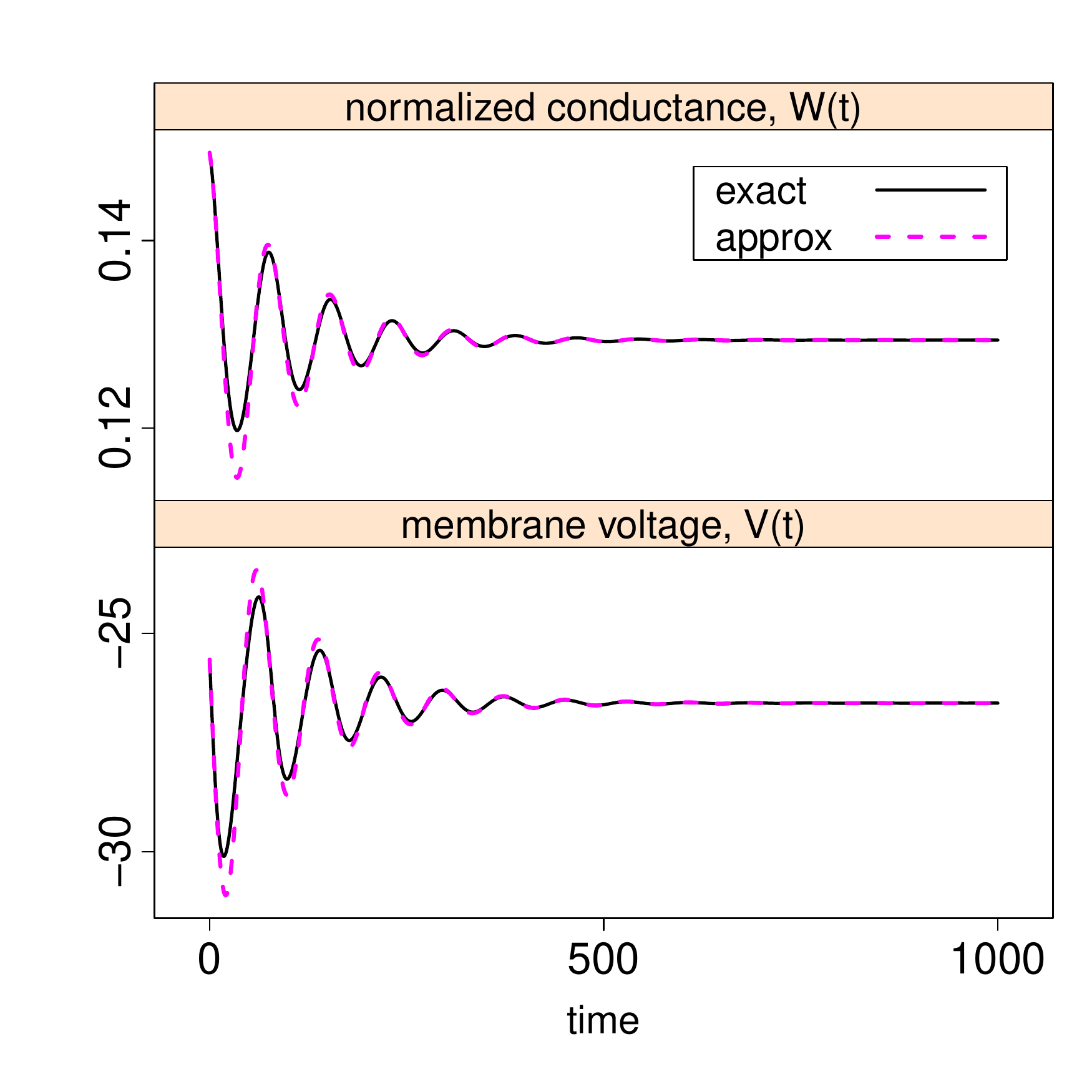}
\caption{ \label{DeterministicApprox} The solution of the deterministic model
\eqref{dV}--\eqref{beta} with $\sigma =0$ (black full drawn curves) is compared
to the linear approximation
\eqref{approxsolution} (cyan dashed curves). Upper panel: normalized conductance $W_t$
 (dimensionless). Lower panel: membrane potential $V_t$
 (mV). Time is measured in ms. \vspace{15mm}}
\end{SCfigure}

\section{Identification of the stochastic process of quiescence}
\label{sec:StocApprox}

In this Section we identify the
stochastic process defined by the linearized system
\eqref{linearsystem} in the limit of small $\lambda$, i.e. under the
condition $\lambda \ll \omega$. The
deterministic system \eqref{approxsolution} has decaying oscillations,
whereas for the stochastic system \eqref{linearsystem}, the noise will
prevent the decay of the oscillations. Can we describe the resulting
process specifically? The answer is that, after a linear change of
variables, this process can be approximated in distribution by a fixed
matrix times a deterministic circular 
motion modulated by an OU process. 

We follow the
development in \cite{BaxendaleGreenwood2010}, where a first step is to
transform the matrix $\mathbf{M}$ into a form which reveals the slow
decay towards the equilibrium point and the fast oscillatory structure
of the deterministic dynamics. Let $\mathbf{Q}$
be a $2 \times 2$ matrix such that
\begin{eqnarray*}
\mathbf{Q}^{-1}\mathbf{M}\mathbf{Q}&=& 
\begin{pmatrix} -\lambda & \omega \\ - \omega & - \lambda 
\end{pmatrix} \, \doteq \,  \mathbf{A}. 
\end{eqnarray*}
A possible choice for $\mathbf{Q}$ is
\begin{eqnarray*}
\mathbf{Q}&=& 
\begin{pmatrix} -\omega & m_{11} + \lambda \\ 0 & m_{21} 
\end{pmatrix} .
\end{eqnarray*}
Let $\tilde{X}_t = \mathbf{Q}^{-1} X_t $, then
\begin{eqnarray}
\label{tildeX}
d\tilde{X}_t &=& \mathbf{A} \tilde{X}_t dt + \mathbf{C}dB_t
\end{eqnarray}
where $\mathbf{C}=\mathbf{Q}^{-1}\mathbf{G}$. 
A further change of variables moves the rotation to form part of the
diffusion coefficient of the linear stochastic system. We define
\begin{eqnarray*}
\tilde{\tilde{X}}_t &=& R_{\omega t} \tilde{X}_t
\end{eqnarray*}
where 
\begin{eqnarray*}
R_{s} &=& \begin{pmatrix}  \cos s &  -\sin s \\   \sin s & \cos s
  \end{pmatrix}
\end{eqnarray*}
is the counterclockwise rotation of angle $s$. Then by Ito's formula
\begin{eqnarray}
\label{dtildetildex}
d\tilde{\tilde{X}}_t &=& -\lambda \tilde{\tilde{X}}_t dt + R_{\omega t} \mathbf{C} dB_t.
\end{eqnarray}
The infinitesimal
covariance matrix in \eqref{tildeX} is 
\begin{eqnarray*}
\mathbf{B} &=& \mathbf{C}\mathbf{C}^T \, = \,
\mathbf{Q}^{-1}\mathbf{G}\mathbf{G}^T (\mathbf{Q}^{-1})^T \, = \, 
\frac{\sigma^2}{m_{21}^2 \omega^2} \begin{pmatrix} (m_{11} + \lambda )^2 &
  \omega (m_{11} + \lambda )\\  
\omega(m_{11} + \lambda ) &  \omega^2 
\end{pmatrix}.
\end{eqnarray*}
Now define
\begin{eqnarray}
\label{tau2}
\tau^2 &=& \frac{1}{2} \mbox{tr} (\mathbf{B}) \, = \,\frac{1}{2}
(B_{11} + B_{22} ) \, = \,  - \frac{\sigma^2 m_{12}}{2 \omega^2 m_{21}},
\end{eqnarray}
where we have used that $ (m_{11} + \lambda )^2 + \omega^2 = -m_{12} m_{21}
$. Finally, we rescale $\tilde{\tilde{X}}_t$ so that we can compare with a
standardized two-dimensional OU process. Let
\begin{eqnarray*}
U_t &=& \frac{\sqrt{\lambda}}{\tau} \, \tilde{\tilde{X}}_{t/\lambda}.
\end{eqnarray*}
Relation \eqref{dtildetildex} becomes
\begin{eqnarray}
\label{dUt}
dU_t &=& -U_t dt + \frac{1}{\tau} R_{\omega t/\lambda} \mathbf{C} d\tilde{B}_t
\end{eqnarray}
where $\tilde{B}_t = \sqrt{\lambda} B_{t/\lambda}$ is another standard
two-dimensional Brownian motion. The following Theorem from
\cite{BaxendaleGreenwood2010} allows us to approximate the process $U_t$
given by \eqref{dUt}, by a two-dimensional OU process with independent
coordinates. 

\begin{thm}
For each fixed $t^*>0$ and $x \in \R^2$ the distribution of $\{ U_t :
0 \leq t \leq t^* \}$ given by \eqref{dUt} with $U_0 = x$ converges
as $\lambda / \omega \rightarrow 0$ to the distribution of the
standardized two-dimensional OU process $\{ S_t : 0 \leq t \leq t^* \}$
generated by
\begin{eqnarray*}
dS_t &=& -S_t dt + dB_t
\end{eqnarray*}
with $S_0 = x$. 
\end{thm}
Here $S_t$ follows a normal distribution, $S_t \sim N \left ( S_0
e^{-t},\frac{1}{2} (1-e^{-2t})  \mathbf{I} \right )$, where $\mathbf{I}$ is the $2\times 2$
identity matrix.
The proof of this Theorem uses a martingale problem convergence
argument and involves the notion of stochastic averaging, where fast
oscillations integrate out revealing the remaining structure
determined by slower oscillations. Another result of this type obtained
by a different method, called multiscale analysis, is in
\cite{KuskeGordilloGreenwood2007}. 

Thus, the process $U_t$ is approximated by $S_t$ if $\lambda \ll \omega$. In our case
$\lambda$ is one order of magnitude smaller than $\omega$. 

Putting together the transformations and the final approximation we
have, in the sense of stochastic process distributions,
\begin{eqnarray}
\label{Xapprox}
X_t &=& \mathbf{Q} \tilde{X}_t \, = \, \mathbf{Q} R_{-\omega t}
\tilde{\tilde{X}}_t  \, = \, 
\mathbf{Q} R_{-\omega t}  \frac{\tau}{\sqrt{\lambda}} \, U_{\lambda t} \, \approx \,
\mathbf{Q} R_{-\omega t}  \frac{\tau}{\sqrt{\lambda}} \, S_{\lambda t}
\nonumber \\
&=& \frac{\tau}{\sqrt{\lambda}} \begin{pmatrix} -\omega & m_{11} + \lambda \\ 0 & m_{21} 
\end{pmatrix} \begin{pmatrix}  \cos
  \omega t &  \sin \omega t \\  - \sin \omega t & \cos \omega
  t   \end{pmatrix}S_{\lambda t}.
\end{eqnarray}
Let us denote by $X_t^a$ the stochastic process on the right hand side
of \eqref{Xapprox}, i.e.  
\begin{eqnarray}
\label{Xa}
X_t^a &=& \tau \mathbf{Q} R_{-\omega t}  \,
S_{\lambda t}/\sqrt{\lambda}.
\end{eqnarray}
To get a sense of how closely the process  $X_t^a$ approximates the
dynamics of the ML process in a neighborhood of
$(V_{\text{eq}},W_{\text{eq}})$ we compare their power spectral
densities, as well as that of the solution of the linearized system
\eqref{linearsystem}. The spectral density of $X_t^a$ and that of$X_t$ satisfying
\eqref{linearsystem} can be calculated explicitly using the power spectrum
formula of \cite{Gardiner} for linear diffusions of the form
\eqref{linearsystem}. In fact $X_t^a$ is 
such a diffusion: the effect of the
stochastic averaging can be seen as replacing $\mathbf{C}$ from
\eqref{tildeX} by a
multiple of the identity in the system \eqref{tildeX}, so the
approximation to $\tilde{X}$ satisfies $d\tilde{X}^a_t = \mathbf{A}
\tilde{X}^a_t dt + \tau dB_t$, where $\tau$ is given by \eqref{tau2}.  
If we transform this equation by $X_t^a = \mathbf{Q}\tilde{X}^a_t $,
we see that $X_t^a$ satisfies
\begin{eqnarray}
\label{XaSDE}
dX^a_t &=& \mathbf{M} X^a_t dt + \tau  \mathbf{Q} dB_t.
\end{eqnarray}
The spectral density of the
first coordinate of $X^a$ is 
\begin{eqnarray*}
S(f) &=& \frac{1}{2\pi} 
\frac{\sigma^2 m_{12}^2}{ \left ( (f^2 -
\mbox{det}(\mathbf{M}))^2 +  (f \mbox{tr}(\mathbf{M}))^2 \right )} \frac{\left ( f^2
+ \mbox{det}(\mathbf{M}) \right )}{2 \omega^2},
\end{eqnarray*}
whereas the spectral density of the
first coordinate of the linearized system, \eqref{linearsystem}, is 
\begin{eqnarray*}
S(f) &=& \frac{1}{2\pi} \frac{\sigma^2 m_{12}^2}{(f^2 - \mbox{det}(\mathbf{M}))^2 + (f \mbox{tr}(\mathbf{M}))^2}.
\end{eqnarray*}
In Fig. \ref{power} the theoretical spectral densities for the two
approximations are plotted, together with the estimated spectral
density of the quiescent process from simulations of the stochastic ML model
\eqref{dV},\eqref{m}--\eqref{beta} and \eqref{dWsigma}. The spectral
density is estimated by 
averaging over at least 20 estimates from paths started at 0 of at
least 450 ms of subthreshold  
fluctuations, and scaled to have the same maximum as the theoretical
spectral density from \eqref{XaSDE}. The averaging is done to reduce
the large variance connected with spectral density estimation,
avoiding any smoothing. Thus, the estimator is approximately
unbiased, see also \cite{MB2005} where this approach is treated. The
estimation is done for $\sigma^* = 0.03, 0.05$ and 
0.1. For higher noise, the lengths of subthreshold fluctuations
between spikes are too short to reliably estimate the spectral
density. Moreover, $\sigma^*=0.1$ corresponds to a number of ion
channels $N \approx 900$, which can be considered a minimum
acceptable number for the diffusion approximation to be relevant. The
value $\sigma^*=0.03$ corresponds to $N \approx 10,000$. Remember that
$\sigma = 0.034 \sigma^*$, see \eqref{G}.

The approximations are only acceptable for small noise, which is
expected, since larger noise brings the process to areas further away
from the fixed point, where non-linearities become increasingly
important. 

\begin{figure}[h]
\centerline{\includegraphics[width=4.5cm]{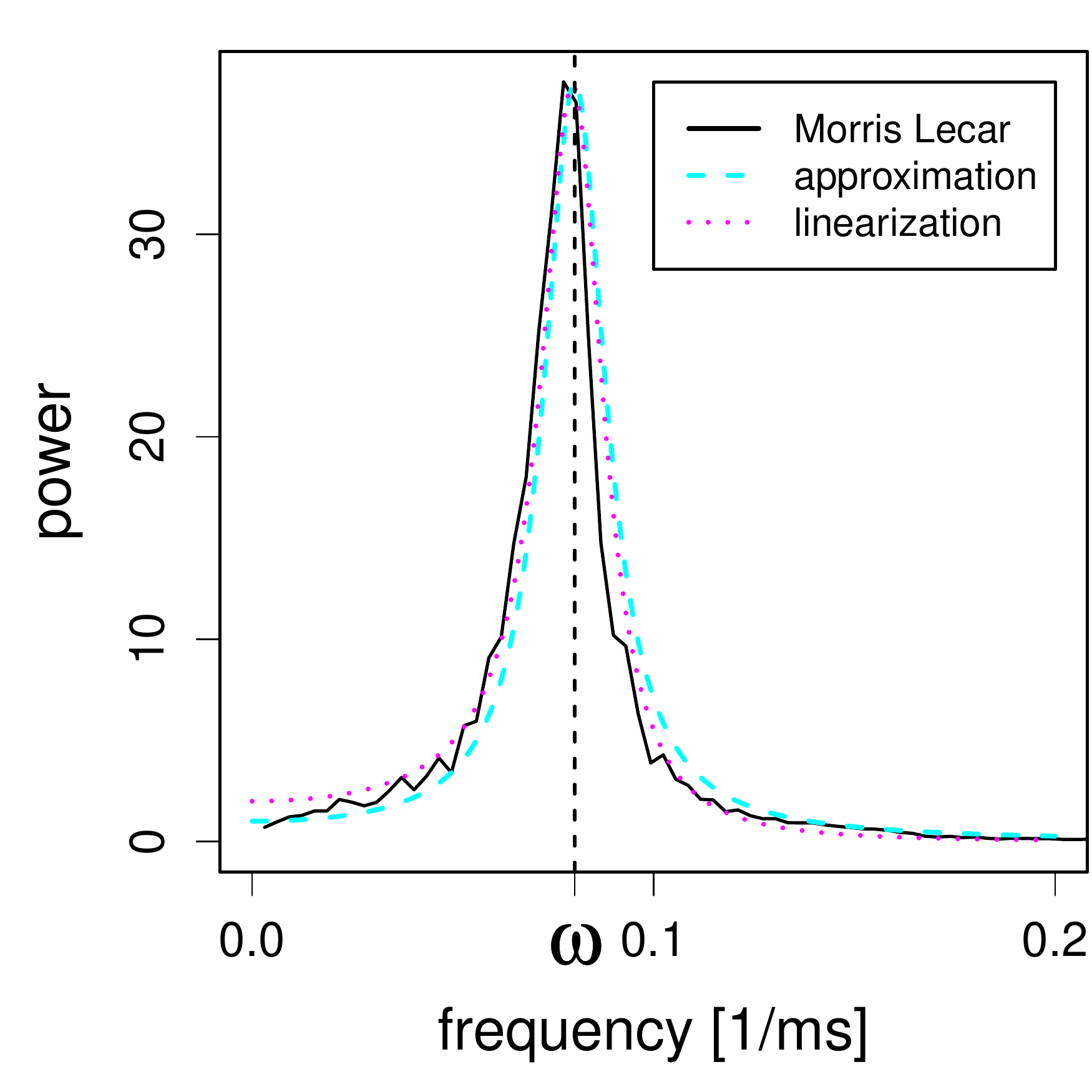}\includegraphics[width=4.5cm]{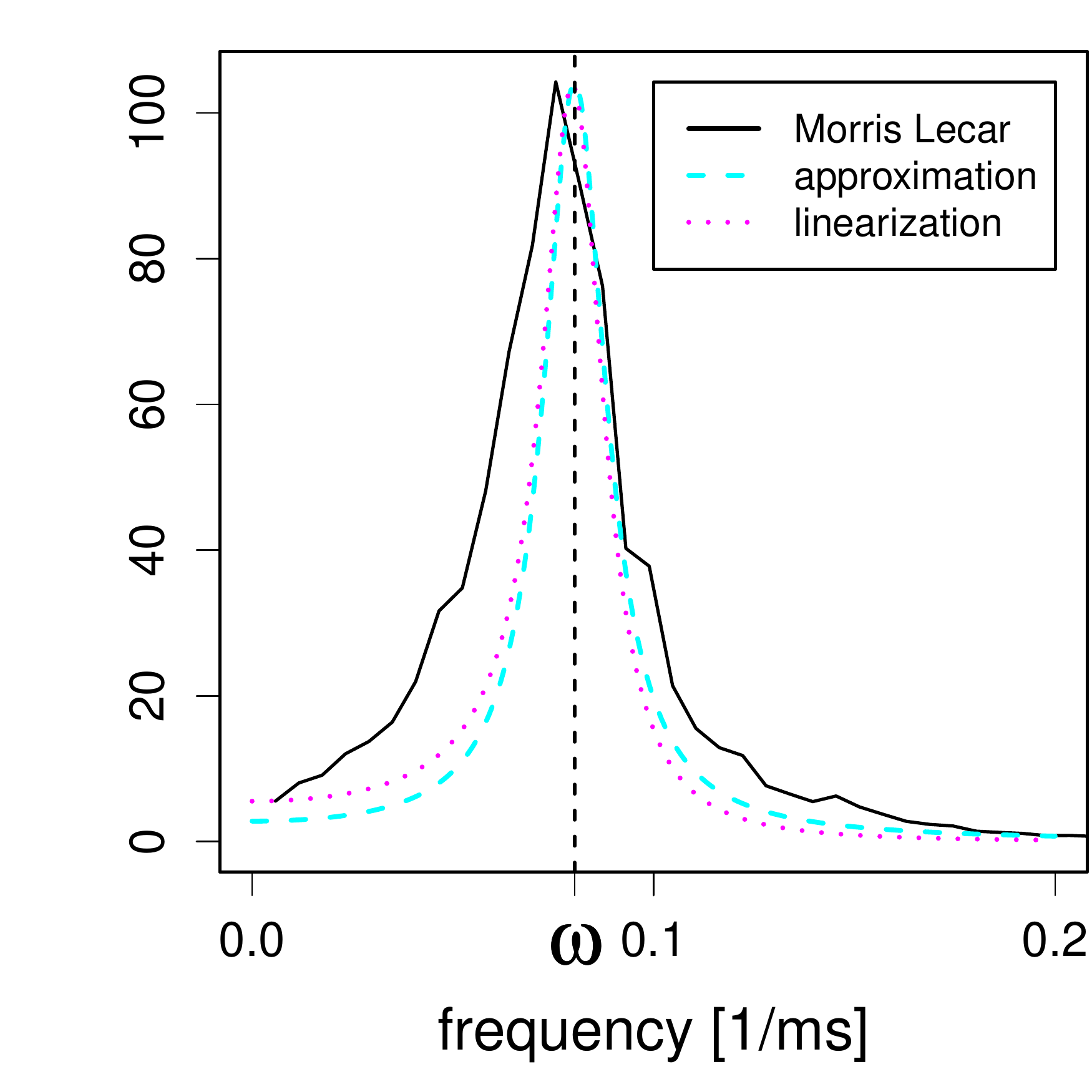}\includegraphics[width=4.5cm]{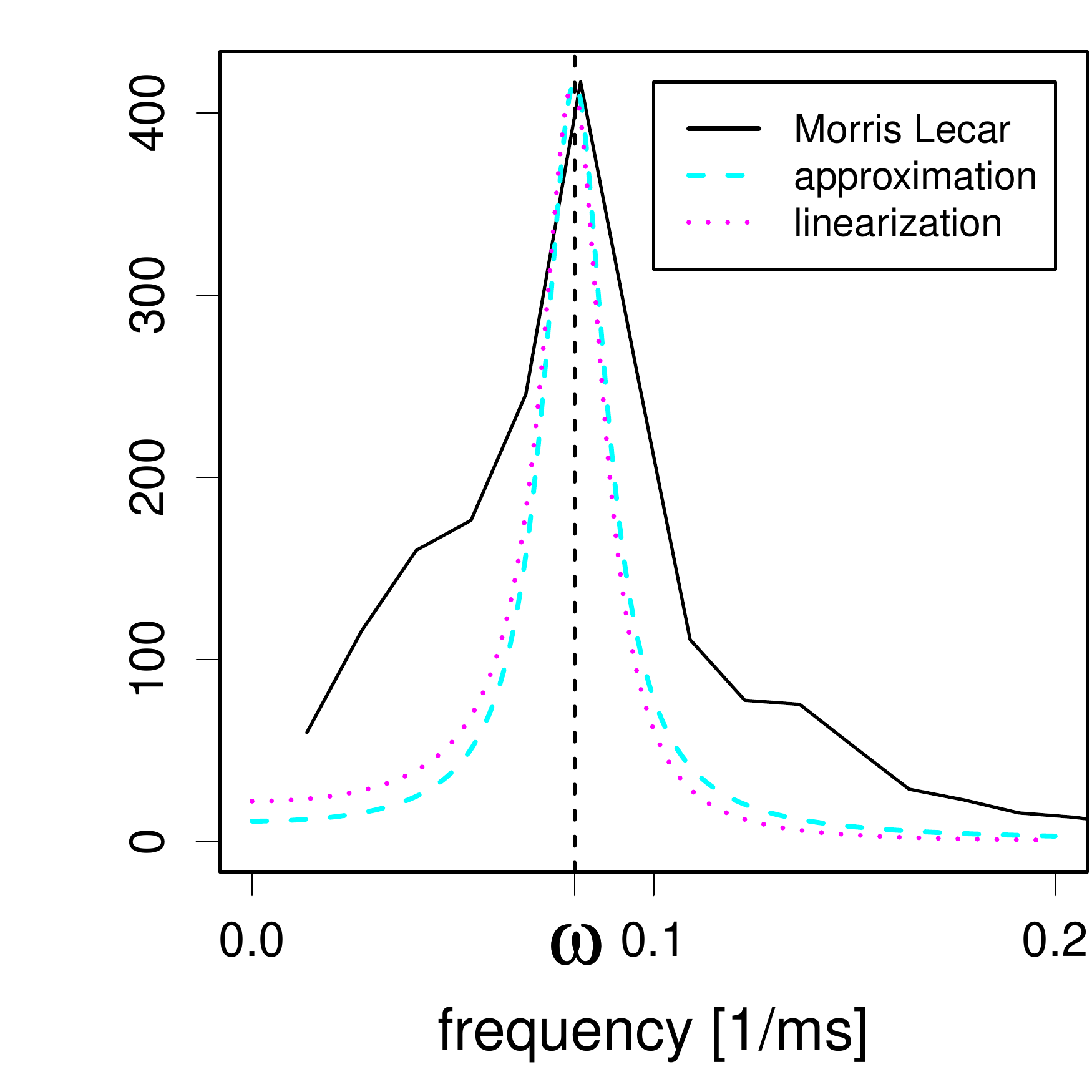}}
\caption{ \label{power} Spectral density estimated from simulations
  between spikes of
  model \eqref{dV}, \eqref{m}--\eqref{beta}, \eqref{dWsigma} (black solid line), theoretical
  spectral density of model \eqref{XaSDE} (cyan dashed line), and  theoretical
  spectral density of model \eqref{linearsystem} (magenta dotted
  line). From left to right: $\sigma^* = 0.03, 0.05$ and $0.1$.   }
\end{figure}

\section{Reconstructing the stochastic ML firing mechanism}
\label{Sec:firing}

In this Section we construct a firing mechanism matching that of the
stochastic ML neuron. In Section \ref{LIFsection} we will define a new
LIF-type process by combining this firing mechanism with the radial OU
process. This new model will, for small $\sigma$, have an ISI
distribution similar to that of the ML.

Firing in model \eqref{dV}, \eqref{m}--\eqref{beta} and
\eqref{dWsigma} occurs when the 
stochastic dynamics shifts from a path circulating the stable
equilibrium, modulated by an OU, to a noisy circuiting of the stable
limit cycle.   
This shift happens, roughly, when the orbit passes from the inside
to the outside of the unstable limit cycle. When the orbit comes close
to the unstable limit cycle, it will follow this limit cycle for a short
time, and then escape either to the inside, i.e. continue
its subthreshold oscillations, or to the outside and a spike will
occur. This understanding is not accurate enough to be implemented
as a firing scheme for the radial OU process \eqref{radialOU}, as we discuss
further in Section \ref{ISIs}. Hence, we embed the process $X^a$
defined by \eqref{Xa} in the stochastic ML model by constructing a
firing mechanism mimicking that of the ML itself. It is clear that in
the ML model, starting inside the unstable limit cycle, a spike will
occur with increasing probability, the further away the process is
from the fixed point. In order to construct a firing mechanism
matching that of ML, we will estimate, from simulations, the
conditional probability that the ML fires, given that the trajectory
of the ML crosses the line $L=\{(v,w) : v=V_{\mbox{eq}},w <
W_{\mbox{eq}} \}$. We computed
estimates from simulated data using crossings of the line $L$ as
follows. 

For a given value of $\sigma^*$ and distance $l$
from the fixed point, a short trajectory  
starting in $(V_{\mbox{eq}},W_{\mbox{eq}}-l)$ was simulated from model
\eqref{dV}, \eqref{m}--\eqref{beta} and \eqref{dWsigma}, and it was
registered whether firing
occurred in the first cycle of the stochastic path around
$(V_{\mbox{eq}},W_{\mbox{eq}}) $. Firing was defined by the path   
crossing the line $v=0$, which is well above the largest level
inside the unstable limit cycle, see
Fig. \ref{DeterministicMorrisLecar}B. This was repeated $1000$
times, and estimates of the conditional probability of spiking,
$\hat{p} (l, 
\sigma^*)$, were computed as the frequency of the
trajectories where firing occurred. The procedure was repeated for
$l=l_i=i\delta, i=1, \ldots, 25$, where 
$\delta$ is the distance to the stable limit cycle divided by 20. In
this way a grid of possible $l$ values was covered, starting from
$l=0$ at the fixed
point, where the probability of firing is close to zero, 
to a point on $L$ below the stable limit cycle, where the probability of firing is
close to one. The estimation was, furthermore, repeated for $\sigma^* =
0.01$ to $0.08$ in steps of $0.01$.

For each fixed $\sigma^*$, the estimates of the conditional 
probability appear 
to depend in a sigmoidal way on the distance from the fixed point. We
assumed the conditional firing probability to be of the form 
\begin{eqnarray}
\label{sigmoidal}
p(l) &=& \frac{1}{1+ \exp((\alpha-l)/\beta)}. 
\end{eqnarray}
\begin{figure}[t]
\centerline{\includegraphics[width=6.5cm]{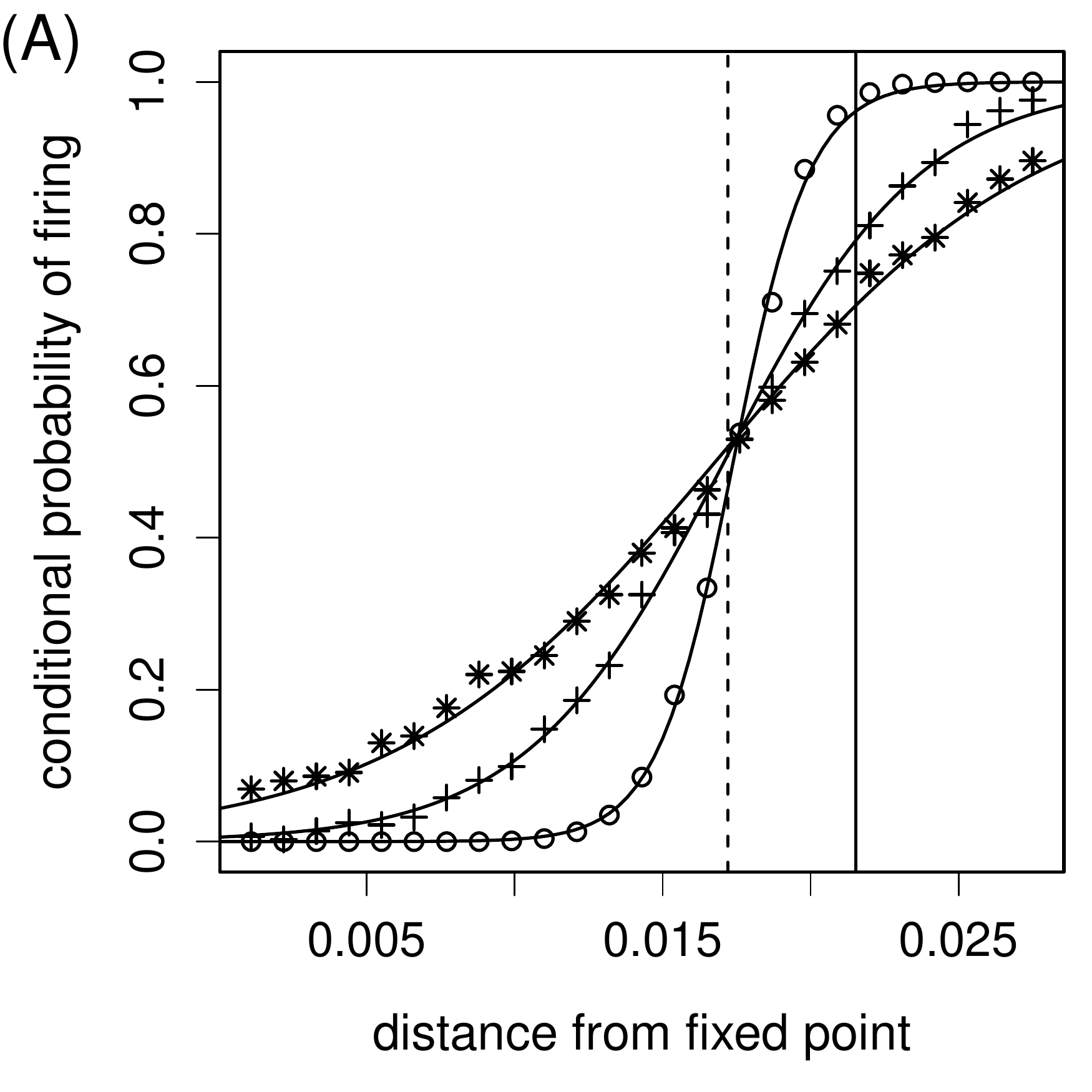}\hspace{2mm}\includegraphics[width=6.5cm]{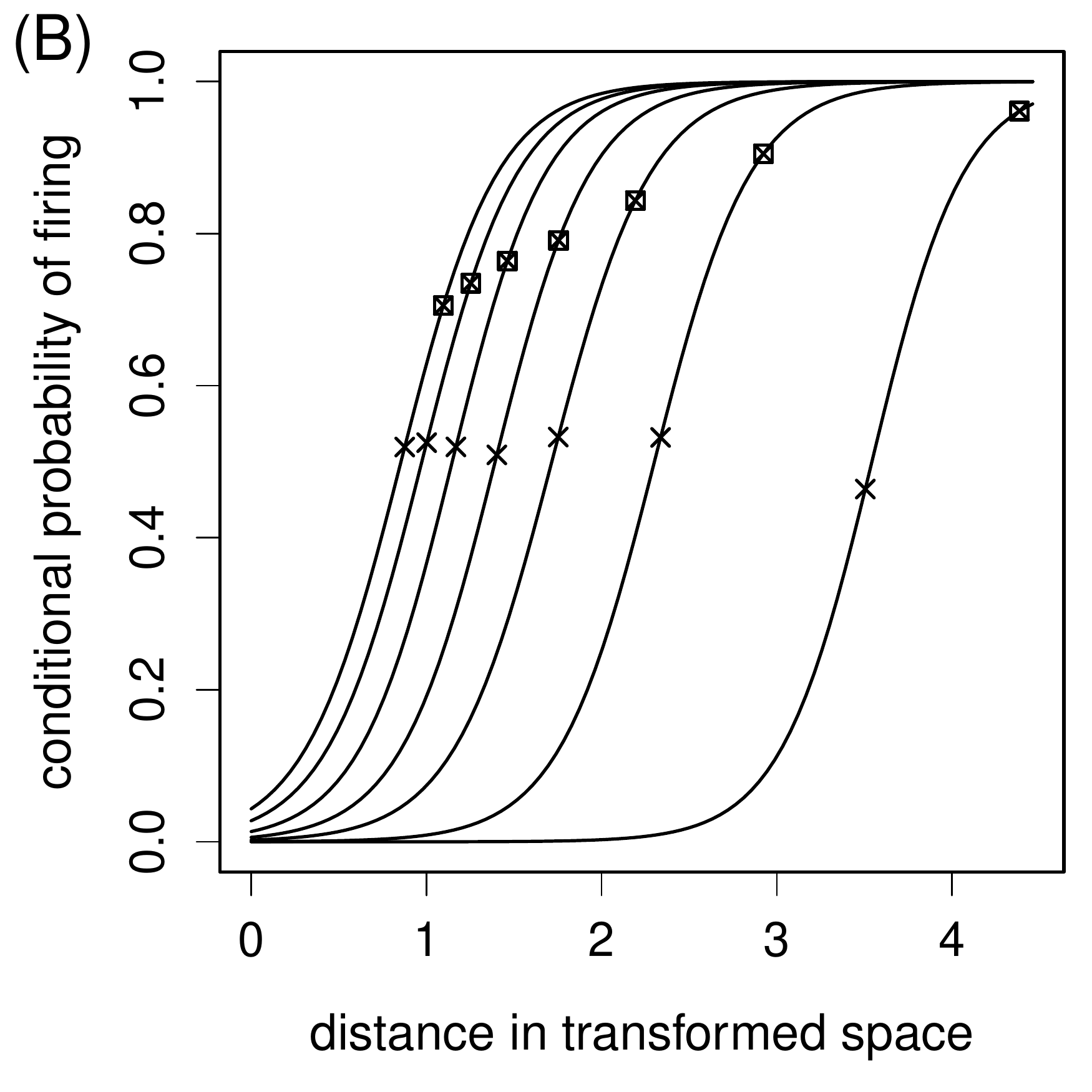}}
\caption{ \label{spikingprobability} Conditional probability of spiking
  when crossing the line $L=\{(v,w) : v=V_{\mbox{eq}},w <
W_{\mbox{eq}} \}$ for different values of $\sigma^*$. (A) Original
space. The circles, plus's and stars are individual 
nonparametric estimates obtained using $\sigma^* = 0.02, 0.05$ and
0.08, respectively, with the fitted curves on top given by
\eqref{sigmoidal}. The dashed line indicates where the unstable limit
cycle 
crosses $L$, the full drawn line where the stable limit cycle crosses
$L$. (B) The fitted curves in the transformed space for $\sigma^* =
0.02, 0.03, 0.04,0.05, 0.06, 0.07$ and $0.08$ (right to left), as a function
of the distance from the fixed point in the transformed
coordinates. The crosses and boxed crosses indicate the crossing 
  of the unstable and stable limit cycles of $L$, respectively, which depend on
  $\sigma = 0.034 \sigma^*$. }
\end{figure}
\begin{table}[b]
  \centering
  \begin{tabular}{c|cccccccc}
$\sigma^*$ &0.01&0.02&0.03&0.04&0.05&0.06&0.07&0.08\\
\hline
$\hat \alpha$ & 0.0174& 0.0174& 0.0169& 0.0168& 0.0171& 0.0169&0.0167&0.0168 \\
$\hat \beta$ & 0.0006 & 0.0013& 0.0020& 0.0028& 0.0033& 0.0039&0.0047&0.0054 \\
$\hat \alpha \sqrt{2 \lambda}/\sigma$ & 7.1022 & 3.5426 & 2.3012 & 1.7156 & 1.3922 & 1.1474 & 0.9739 & 0.8549\\
$\hat \beta \sqrt{2 \lambda}/\sigma$ & 0.2590 & 0.2624 & 0.2759 & 0.2831 & 0.2718 & 0.2674 & 0.2738 & 0.2764 \\
\hline
  \end{tabular}
  \caption{Estimates of regression parameters for $p(\cdot)$ in the
    original space (first two rows), and in the transformed
    coordinates (last two rows).}
  \label{tab:regressionestimates}
\end{table}
The parameters $\alpha$ and $\beta$ were estimated using non-linear
regression of the 25 estimates of $\hat 
p(l_i;\sigma^*)$ on $l$.
In Fig. \ref{spikingprobability}A these parametric estimates are
plotted, as well as the
individual 
nonparametric estimates $\hat p$ for $\sigma^* = 0.02, 0.05$ and
$0.08$.  We see that the family of estimates, $\hat{p}$, fits the
hypothetic 
curve quite well for each value of $\sigma^*$. 
Regression estimates are reported in Table 
\ref{tab:regressionestimates}. 
Note that $\alpha$ is the distance along $L$ from $W_{\mbox{eq}}$ 
at which the conditional probability of firing equals one half. For all
values of $\sigma^*$, the estimate of $\alpha$ is close to the
distance along $L$ between $W_{\mbox{eq}}$ and  
the unstable limit
cycle, which equals 0.0172. In other words, the probability of firing,
if the path starts at the intersection of $L$ with the unstable limit
cycle, is about 1/2.  
The parameter $\beta$ indicates the width of a band around $\alpha$
where the conditional probability essentially changes. For instance,
if $l \in 
\alpha \pm \beta$ then $p(l) \in (0.27,0.73)$, if $l \in
\alpha \pm 2\beta$ then $p(l) \in (0.12,0.88)$. As expected, the
estimate of $\beta$ increases with increasing $\sigma^*$, and for
small noise the conditional probability approaches a step function
since the process is mostly dominated by the drift. A step function would
correspond to the firing being represented by a first-passage time of
a fixed threshold. Note though that $\hat \beta$ is approximately  
proportional to $\sigma^*$, and thus, as we said earlier and will see
in the following, 
a fixed threshold at the crossing of the unstable limit cycle does not
reproduce the desired spiking characteristics.

In order to simplify the construction in Section \ref{LIFsection} of a
LIF model 
which, together with a firing rule, behaves like the stochastic ML, we
will change coordinates as follows. 
Observe that \eqref{Xa} can be written
\begin{eqnarray}
\label{QXapprox}
\frac{\sqrt{\lambda}}{\tau} \mathbf{Q}^{-1} X_t^a  &=& \begin{pmatrix}  \cos
  \omega t &  \sin \omega t \\  - \sin \omega t & \cos \omega
  t   \end{pmatrix} S_{\lambda t},
\end{eqnarray}
so for fixed $t$, $\sqrt{\lambda}
\mathbf{Q}^{-1} X_t^a/\tau $ is the clockwise 
rotation by angle $\omega t$ of the orthogonal pair $(S_{\lambda
  t}^{(1)},S_{\lambda t}^{(2)})$. We define a transformation of the
space $(v,w)$ by centering at $(V_{\mbox{eq}},W_{\mbox{eq}})$ and
normalizing as in \eqref{QXapprox}. Let 
\begin{eqnarray}
\label{transformedspace}
\begin{pmatrix}  \tilde v \\ \tilde w   \end{pmatrix} &=&
\frac{\sqrt{\lambda}}{\tau} \mathbf{Q}^{-1} \begin{pmatrix}
  v-V_{\text{eq}} \\ w-W_{\text{eq}}   \end{pmatrix}  
\end{eqnarray}
be the coordinates of the transformed space. In the new coordinates
our process is simplified to a rotation modulated by a standard
two-dimensional OU process with independent components. 

The transformation depends on $\sigma = 0.034 \sigma^*$, namely, the transformed unstable
limit cycle becomes smaller with increasing noise, through the value
of $\tau$ given in \eqref{tau2}. This is exactly
what is causing a higher firing probability for 
larger $\sigma^*$. The line $L$ will in the transformed
space be
\begin{eqnarray*}
\tilde L &=&
\frac{\sqrt{\lambda}}{\tau} \mathbf{Q}^{-1} \begin{pmatrix}
  0 \\ l \end{pmatrix} =  \frac{\sqrt{\lambda}}{m_{21}\tau} \begin{pmatrix}
  \frac{m_{11}+\lambda}{\omega} \\ 1 \end{pmatrix} l 
\end{eqnarray*}
for $l \geq 0$. A distance $l$ will thus transform to a distance $r =
(\sqrt{2\lambda}/\sigma) l$, and the conditional probability of
firing \eqref{sigmoidal} transforms to
\begin{eqnarray}
\label{sigmoidaltransformed}
p(r) &=& \frac{1}{1+
  \exp((\alpha^*-r)/\beta^*)},  
\end{eqnarray}
where
$\alpha^* = \alpha \sqrt{2\lambda}/\sigma$ and
$\beta^*=\beta\sqrt{2\lambda}/\sigma$.  
The fitted curves of \eqref{sigmoidaltransformed} for $\sigma^* = 0.02
- 0.08$, as a function 
of the distance from the fixed point in the transformed
coordinates are given in Fig. \ref{spikingprobability}B, with
indication of the crossings 
  of the unstable and stable limit cycles, respectively, which now depend on
  $\sigma$. Note that in the transformed space, the width of the
  band where the conditional probability is essentially different from
  0 or 1 is nearly constant, see Table
  \ref{tab:regressionestimates}. From here on we use the coordinates
  defined by \eqref{transformedspace}. 

\section{Construction of a leaky-integrate-and-fire model with ML
  firing statistics}
\label{LIFsection}

The simpler stochastic LIF models sacrifice realism for mathematical tractability
\cite{Burkitt2006,GerstnerKistler2002}. In these models, a neuron is
characterized by a single stochastic differential equation
describing the evolution of neuronal membrane potential 
depending on time, 
\begin{eqnarray}
\label{SDEmodel} 
d X_t &=& \mu (X_t ) d t + \sigma (X_t ) d B_t ,\quad\quad
X_0 =  x_0,
\end{eqnarray}
where $X_t$ corresponds to
$V_t$ in the ML model, together with a threshold firing rule,
\begin{eqnarray}
\label{TS}
T&=&\inf\{t>0:X_t\geq S \}.
\end{eqnarray}
In this Section we define a LIF model which does not make this
compromise, using the result of Section \ref{sec:StocApprox} and the
firing mechanism defined in Section \ref{Sec:firing}.

The distance of the approximate 
process $\sqrt{\lambda}
\mathbf{Q}^{-1} X_t^a/\tau$ of \eqref{QXapprox} from the point $(0,0)$ at
time $t$ is 
given by the modulus of the two-dimensional 
standardized OU process $S_{\lambda t}$. 
The modulus of $S_{\lambda t}$ at time $t$ is given by the process
\begin{eqnarray*}
R_{\lambda t} &=& \sqrt{(S_{\lambda t}^{(1)})^2 + (S_{\lambda t}^{(2)})^2},
\end{eqnarray*}
which is a standard radial OU process with two degrees of freedom. It
has state space $(0,\infty)$, and solves the stochastic differential
equation 
\begin{eqnarray}
\label{radialOU}
dR_{\lambda t} &=& \left ( \frac{1}{2R_{\lambda t}} - R_{\lambda t} \right ) dt + dW_{\lambda t}, 
\end{eqnarray}
see e.g. \cite{BorodinSalminen2002}. We define a new LIF process by
\eqref{radialOU}, and firing 
mechanism derived from \eqref{sigmoidaltransformed}. After each firing, we
will reset the time to 0 and assume the process reset to 0, i.e. $R_0
= 0$, corresponding to 
$S_0 = (0,0)$ and $(V_0,W_0)=(V_{\mbox{eq}},W_{\mbox{eq}})$.  
By Ito's formula, the process $Y_u = R_u^2$ satisfies 
the stochastic differential equation 
\begin{eqnarray}
\label{Yu}
dY_u &=& 2\left ( 1 - Y_u \right ) du + 2\sqrt{Y_u}dW_u, 
\end{eqnarray}
and is thus a square-root process, see e.g. \cite{CIR85}, also called
a Feller or a 
Cox-Ingersoll-Ross process. This process is ergodic, and its stationary
distribution is the exponential distribution with mean one. It follows
that the stationary distribution of $R_u$ has density $f(r) =
2re^{-r^2}$ on $(0, \infty)$, i.e. it follows a Rayleigh distribution. The
transition density of $Y_u$ starting at $y_0$ at time 0, is a
non-central $\chi^2$-distribution with two 
degrees of freedom and non-centrality parameter $\delta (u, y_0) =
2y_0e^{-2u}/(1-e^{-2u})$. Then $2Y_u/(1-e^{-2u})$ follows the
standard non-central $\chi^2$-distribution $F_{\chi^2}
(2y/(1-e^{-2u}),2,\delta (u, y_0))$. It is particularly simple because
of the integer degrees of freedom. Transforming to the radial OU
we obtain the transition density of $R_u$ starting at $s$ at time 0
\begin{eqnarray}
\label{transitiondensityR} 
f_u(r,s) &=& \frac{2r}{1-e^{-2u}} \exp \left \{
  -\frac{r^2+s^2e^{-2u}}{1-e^{-2u}} \right \} I_0 \left (
  \frac{rs}{\sinh (u)} \right ),
\end{eqnarray}
where $I_0 (x)= \frac{1}{\pi} \int_0^{\pi} e^{x \cos \theta} d \theta$
is the modified Bessel function of the first kind of 
index 0.

Writing the two-dimensional process $S_u$ in polar coordinates, $R_u$ and
$\theta_u$, where  $\theta_u$ is the angle at time $u$ to the positive
part of the first coordinate, we find that the modulus and the angle
are independent, and that $\theta_u$ is uniformly distributed on
$(0,2\pi)$. This can e.g. be seen from the fact that $S_u^{(1)}$ and
$S_u^{(2)}$ are independent normal with mean 0 and equal 
variances. Thus, for fixed $u$, $S_u^{(2)}/S_u^{(1)}$ is standard
Cauchy distributed and $\theta_u= \arctan (S_u^{(2)}/S_u^{(1)})$ is
$U(0,2\pi)$. 

Let $T$ denote the firing time random variable. 
We want to compute the density
of the distribution of $T$, and for this we find it convenient to
express this density in terms of the conditional hazard rate, 
\begin{eqnarray*}
\alpha(t,r) &=& \lim_{\Delta t \rightarrow 0}\frac{1}{\Delta t} P(t
\leq T < t + \Delta t \, | \,  T \geq t, R_{\lambda t} = r).
\end{eqnarray*}
This function is the density of the conditional probability, given the
position on $L$ is $r$ at time $t$,  of a spike occuring in the next
small time interval, given that it has not yet occurred. 

From standard results from survival analysis, see
e.g. \cite{Aalenbook}, we obtain
\begin{eqnarray*}
P(T> t \, | \, R_{\lambda s}, 0\leq s \leq t)&=& \exp \left ( -\int_0^t
  \alpha(R_{\lambda s})ds \right ) .
\end{eqnarray*}
The unconditional distribution is then given by 
\begin{eqnarray}
\label{St} 
P(T> t )&=& E \left (\exp \left ( -\int_0^t
  \alpha(R_{\lambda s})ds \right ) \right )
\end{eqnarray}
where $E(\cdot)$ denotes expectation with respect to the distribution
of $R$. The density is thus
\begin{eqnarray}
\label{gt} 
g(t) \, = \, \frac{d}{dt} P(T\leq  t )&=& E \left ( \alpha (R_{\lambda
    t}) \exp 
  \left ( -\int_0^t 
  \alpha(R_{\lambda s})ds \right ) \right ).
\end{eqnarray}
The firing is defined to be initiated from $L$, and on average the
process crosses $L$ every $2\pi/\omega = 78.2$ time units. Using
\eqref{sigmoidaltransformed}, the estimated conditional probability of
firing given the position on $L$ is $r$, which by definition does not
depend on $t$, we estimate the hazard 
rate as 
\begin{eqnarray}
\label{sigmoidaltransformedalpha}
\alpha (t,r) \, = \, \alpha (r) &=& \frac{\omega}{2\pi} \frac{1}{1+
  \exp((\alpha^*-r)/\beta^*)}.  
\end{eqnarray}
Note that it is bounded. This is not realistic, since a very large
value of $r$ should cause immediate firing. 

In \cite{ISI:000237628500003} a firing rule with unbounded hazard rate
was proposed, and in \cite{Jahn2010} it was shown to fit well to
experimental data. Therefore, 
we will also see how our model performs if we use in the
firing mechanism a hazard rate of the form
\begin{eqnarray}
\label{expalpha} 
\alpha(t,r) \, = \, \alpha(r) &=& \exp((r-\alpha)/\beta)
\end{eqnarray}
for $\alpha, \beta >0$. Like before, $\alpha$ plays the role of a
threshold, 
and $\beta$ gives the 
width of the threshold region. When $\beta \rightarrow 0$, the firing
rule converges to a fixed threshold crossing. To estimate
$\alpha$ and $\beta$ in \eqref{expalpha}, we simulated 1000 spike
times from the ML. The cumulative hazard $A(t) = \int_0^{t} \alpha(t)$ 
was then estimated from the simulated spike times by the standard
empirical Nelson-Aalen 
estimator. The theoretical cumulative hazard using \eqref{radialOU} and
\eqref{expalpha} can be calculated as
\begin{eqnarray}
\label{At} 
A(t) &=& E \left ( \int_0^{t} \alpha(R_{\lambda s}) ds \right ) \, =
\, \exp \left (-\frac{\alpha}{\beta} \right ) \int_0^{t} E \left
  (\exp \left (\frac{R_{\lambda s} }{\beta} \right ) \right ) ds
\nonumber \\
&=& \sqrt{\pi} \exp \left (-\frac{\alpha}{\beta} \right ) \int_0^{t}
\left ( g(s)   \exp \left (
    \frac{1}{4} g(s)^2\right ) \Phi \left (
      g(s) \right ) +1\right 
) ds  
\end{eqnarray}
where we have used the density $f_{\lambda s}(r,0)$ given in
\eqref{transitiondensityR}. Here, $g(s) = \sqrt{1-e^{-2\lambda
    s}}/\beta$, and $\Phi (\cdot)$ is the standard 
normal cumulative distribution function. Then, $\alpha$ and $\beta$
were estimated by the least square distance between \eqref{At} and the
estimated cumulative hazard from the simulated spike times.
For $\sigma^*=0.05$ the estimates 
were $\alpha = 6.31$ and $\beta = 0.76$. 

The final model is
\begin{eqnarray}
\label{radialOUwithpoisson}
dR_{u} &=& \left ( \frac{1}{2R_{u}} - R_{u} \right ) dt + dW_{u} -
R_{u-} \mu (R_{u-},du), 
\end{eqnarray}
where $\mu (R_{u-},du)$ is a Poisson measure with intensity
$\alpha(R_{u-})$, and $R_{u-}$ denotes the left limit of $R_u$. Here,
$\alpha (\cdot)$ is either given by \eqref{sigmoidaltransformedalpha}
or \eqref{expalpha}. The 
jump size is $-R_{u-}$, thus giving the reset to 0 at spike times. 

A reasonable alternative to the soft threshold firing mechanism used
here would be to use the firing rule defined by a threshold as in
the classical LIF models, equation \eqref{TS}. A natural choice of
threshold would be where the 
LIF process reaches a level corresponding to the unstable limit
cycle. In fact, according to 
our estimates in Fig. \ref{spikingprobability} and Table
\ref{tab:regressionestimates}, the 
firing probability 
of the ML at this threshold is around $1/2$. However, the ISI
distribution estimated from simulations 
using a hard threshold at the unstable limit cycle is shifted towards
larger times, relative to the ML ISI distribution. This happens
because the process might cycle many times inside the unstable limit
cycle, so even if the probability of spiking in a single cycle is
small, the total probability is not negligible. This is lost
when only a hard threshold is considered. Instead we chose the
threshold value such that the mean of the ML ISI distribution and the
mean of the LIF ISI distribution were the same. In
\cite{ISI:000255902700007}, the mean of $T$ from \eqref{TS} with $X_t = R_t$ 
started at $R_0=0$ is given using a hypergeometric function,
\begin{eqnarray}
\label{ET}
E(T) &=& \frac{S^2}{2}  \,  _2F_2 \left (1,1;2,2;S^2 \right ).
\end{eqnarray}
The average of the 1000 ML firing times for $\sigma^*=0.05$ was
447. Equating with \eqref{ET} gives a value $S=2.97$ for the hard
threshold. Note that this is much smaller than the estimated $\alpha$
from \eqref{expalpha}.

\section{Comparison of firing statistics}
\label{ISIs}

One of the major issues in computational neuroscience is to determine
the ISI distribution. We therefore simulated the ML model given by \eqref{dV},
  \eqref{m}--\eqref{beta} and \eqref{dWsigma} until spiking, and
  thereafter reset to
  the fixed point. This was done 1000 times, and the time of the
  firing was recorded. The ISI distribution from our
  approximate model is given by the density \eqref{gt}, or
  equivalently, from the survival function \eqref{St}. Due to the law
  of large numbers and since we know
  the exact distribution of $R_u$, for fixed $t$ we can numerically
  determine \eqref{gt} up to any desired precision by choosing $n$ and
  $M$ large enough through the expression
\begin{eqnarray}
\label{gtapprox} 
g(t) &\approx& \frac{1}{M} \sum_{m=1}^M  \alpha \left (
    R_{\lambda t}^{(m)} \right ) \exp
  \left ( - \frac{t}{n} \sum_{i=1}^n 
  \frac{\alpha \left ( R_{i\lambda t/n}^{(m)} \right ) +\alpha \left (
      R_{(i-1)\lambda t/n}^{(m)} \right )  }{2} \right ) . 
\end{eqnarray}
Here $(R_{0}^{(m)}, \ldots , R_{i\lambda t/n}^{(m)}, \ldots ,
R_{\lambda t}^{(m)})$
are $M$ realizations of $R_{i\lambda t/n}$, $i=0,1, \ldots , n$, and
the integral has 
been approximated by the trapezoidal rule. The hazard rate is either
given by \eqref{sigmoidaltransformedalpha} or \eqref{expalpha}.

 The results are illustrated in Figure
  \ref{ISIdensities} for $\sigma^*=0.05$, using $M=1000$. The
  estimated ISI 
  distributions from our 
  approximate model with both firing mechanisms compare well with the
  estimated ISI distribution 
  of ML reset to 0 after firings. On the contrary, the hard threshold
  does not reproduce the ISI distribution well, e.g. the right tail is
  too heavy. This is because the probability of firing during low
  subthreshold activity is set to 0, whereas we have seen it is
  not. 

\begin{SCfigure}
  \centering
\includegraphics[width=8cm]{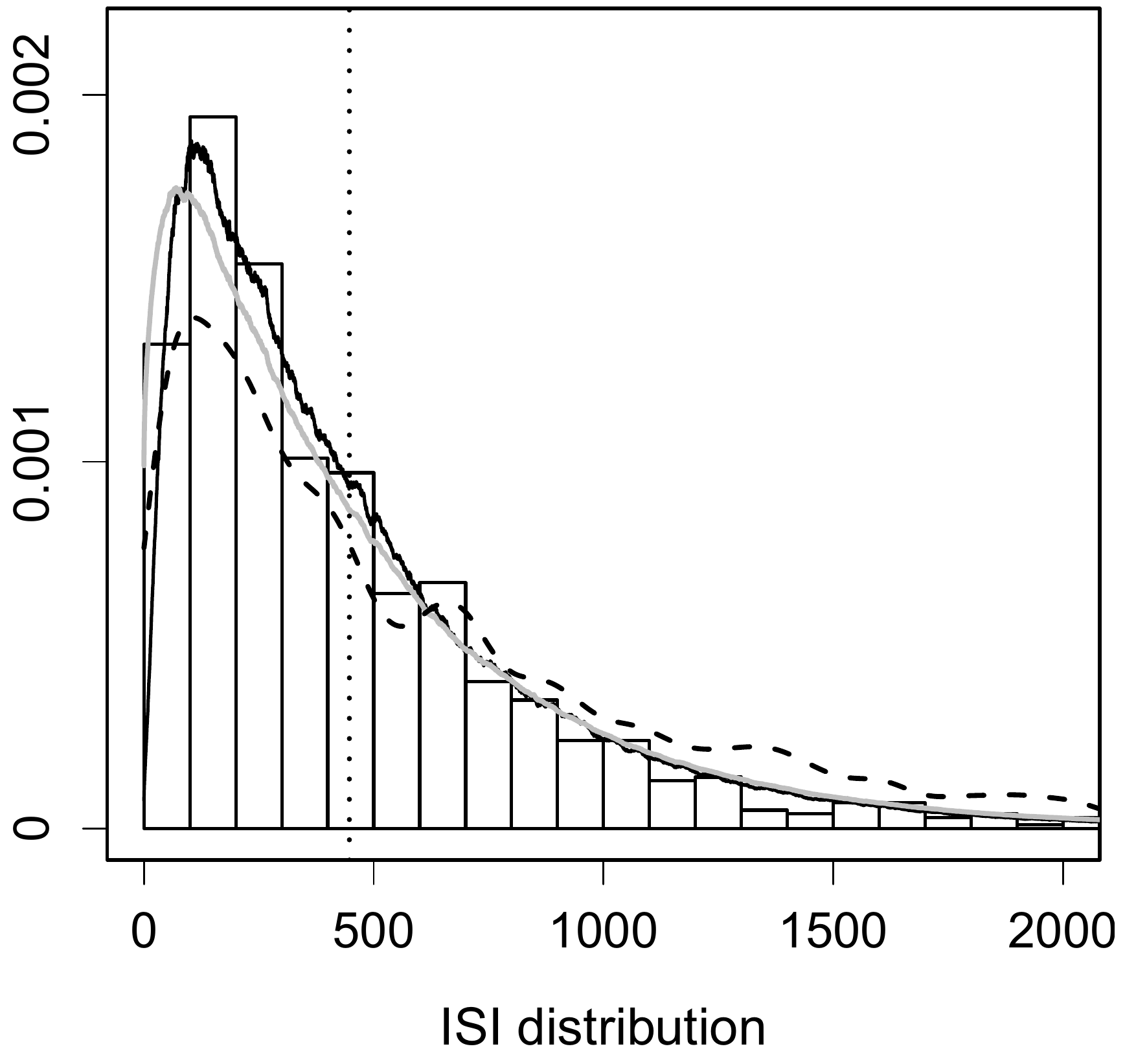}
\caption{ \label{ISIdensities} Distribution of firing times for
  $\sigma^* = 0.05$. The histogram is  
  based on 1000 simulated firing times from the ML model, the vertical dotted
  line is the average. Curves are
  estimates of the probability density, equation
  \eqref{gtapprox}. Black curve 
  is estimated using \eqref{sigmoidaltransformedalpha}, gray curve is
  estimated using \eqref{expalpha}, dashed curve is estimated using a
  fixed threshold, \eqref{TS}.
  \vspace{3mm} }
\end{SCfigure}

\section{Discussion}
\label{discussion}

A stochastic LIF model constructed with a radial OU process and firing
mechanism of either logistic or exponential type has been shown to
mimic the ISI statistics 
of a ML neuron model. It captures
subthreshold dynamics, not of the membrane potential alone, but of a
combination of the membrane potential and ion channels. This construction will allow us to answer
several questions about ML models, which have been accessible only for
LIF models, even though the latter have less biological motivation.

An example of such a question would be: Using ISI experimental data, the
noise standard deviation $\sigma$ can be estimated
\cite{LanskyDitlevsen2008}. In principle, this should also be possible
from our new LIF model, even though we use a soft threshold. This will
give an estimate of $N$, the number of ion channels involved,
through the relation $(\sigma^*)^2 \approx 9/N$.

A question we have not explored is: what is the best way to
restart our new LIF model? In our simulations we restarted both our
LIF and the ML at the fixed point of the ML. However, an uninterrupted
stochastic ML produces continuous paths as in
Fig. \ref{DeterministicMorrisLecar}B. After firing, which means
traversing the large stable limit cycle, possibly several times, they
reenter a neighborhood of the fixed point from its edge. A further
refinement of our LIF model will be obtained by introducing a reentry
mechanism, which mimics this aspect of the ML.


\appendix

\section{Linearization matrix}

The expression for $\mathbf{M}$ in \eqref{linearsystem} is
\begin{eqnarray*}
\mathbf{M} &=&   \begin{pmatrix} m_{11} &
  -g_kW_{\text{eq}}(V_{\text{eq}}-V_K) /C \\[2mm]
2V_{\text{eq}}W_{\text{eq}}\beta \left ( V_{\text{eq}}  \right )/V_4  &
  - \alpha \left ( V_{\text{eq}}  \right )  \end{pmatrix},
\\[2mm]
m_{11} &=& - \frac{V_{\text{eq}}}{C} \left (
\frac{2g_{\text{Ca}}(V_{\text{eq}}-V_{\text{Ca}})\alpha \left (
  V_{\text{eq}}  \right )\beta \left ( V_{\text{eq}}  \right )}{V_2 (\alpha \left (
  V_{\text{eq}}  \right ) + \beta \left ( V_{\text{eq}}  \right ) )^2}
   +  g_{\text{Ca}}m_{\infty}(V_{\text{eq}}) + g_K
  W_{\text{eq}} + g_L \right )
\end{eqnarray*}



\acks
S. Ditlevsen supported by the
Danish Council for Independent Research$\, | \, $Natural
Sciences. P. Greenwood supported by the Statistical and Applied
Mathematical Sciences 
Institute, Research Triangle Park, N.C., and the Mathematical,
Computational and Modeling Sciences Center at Arizona State
University. The Villum Kann Rasmussen foundation supported a 4 months
visiting professorship for P. Greenwood at University of Copenhagen.

%
%
%
%

\bibliographystyle{apt}
\bibliography{bibliography}

\begin{thebibliography}{10}

\bibitem{Aalenbook}
{\sc Aalen, O.~O., Borgan, {\O}. and Gjessing, H.~K.} (2010).
\newblock {\em Survival and Event History Analysis. A process point of view}.
\newblock Springer, New York.

\bibitem{BaxendaleGreenwood2010}
{\sc Baxendale, P. and Greenwood, P.} (2011).
\newblock Sustained oscillations for density dependent {M}arkov processes.
\newblock {\em J. Math. Biol.\/} To appear, available online.

\bibitem{BorodinSalminen2002}
{\sc Borodin, A.~N. and Salminen, P.} (2002).
\newblock {\em Handbook of {B}rownian motion - {F}acts and {F}ormulae}.
\newblock Probability and its applications. Birkhauser Verlag, Basel.

\bibitem{Burkitt2006}
{\sc Burkitt, A.~N.} (2006).
\newblock A review of the integrate-and-fire neuron model: I. homogeneous
  synaptic input.
\newblock {\em Biol. Cybern.\/} {\bf 95,} 1--19.

\bibitem{CIR85}
{\sc Cox, J.~C., Ingersoll, J.~E. and Ross, S.~A.} (1985).
\newblock A theory of the term structure of interest rates.
\newblock {\em Econometrica\/} {\bf 53,} 385--407.

\bibitem{DayanAbbott2001}
{\sc Dayan, P. and Abbott, L.~F.} (2001).
\newblock {\em Theoretical Neuroscience: Computational and Mathematical
  Modeling of Neural Systems}.
\newblock MIT Press, Cambridge.

\bibitem{DitlevsenJacobsen2011}
{\sc Ditlevsen, S. and Jacobsen, M.} (2011).
\newblock Ergodic solutions to multidimensional diffusions.
\newblock In preparation.

\bibitem{MB2005}
{\sc Ditlevsen, S., Yip, K.-P. and Holstein-Rathlou, N.-H.} (2005).
\newblock Parameter estimation in a stochastic model of the tubuloglomerular
  feedback mechanism in a rat nephron.
\newblock {\em Mathematical Biosciences\/} {\bf 194,} 49--69.

\bibitem{PearsonDiffusions}
{\sc Forman, J.~L. and S{\o}rensen, M.} (2008).
\newblock The {P}earson diffusions: A class of statistically tractable
  diffusion processes.
\newblock {\em Scand. J. Stat.\/} {\bf 35,} 438--465.

\bibitem{Gardiner}
{\sc Gardiner, C.~W.} (1990).
\newblock {\em Handbook of stochastic methods for physics, chemistry and the
  natural sciences} 2nd~ed.
\newblock Springer, Berlin Heidelberg.

\bibitem{GerstnerKistler2002}
{\sc Gerstner, W. and Kistler, W.~M.} (2002).
\newblock {\em Spiking Neuron Models}.
\newblock Cambridge Uni. Press, Cambridge.

\bibitem{ISI:000255902700007}
{\sc Graczyk, P. and Jakubowski, T.} ({2008}).
\newblock {Exit times and Poisson kernels of the Ornstein-Uhlenbeck diffusion}.
\newblock {\em {Stochastic Models}\/} {\bf {24},} {314--337}.

\bibitem{HH1952}
{\sc Hodgkin, A.~L. and Huxley, A.~F.} (1952).
\newblock A quantitative description of ion currents and its applications to
  conduction and excitation in nerve membranes.
\newblock {\em J. Physiol.\/} {\bf 117,} 500--544.

\bibitem{BookIz}
{\sc Izhikevich, E.~M.} (2007).
\newblock {\em Dynamical systems in neuroscience: the geometry of excitability
  and bursting}.
\newblock MIT Press, Cambridge.

\bibitem{Jahn2010}
{\sc Jahn, P., Berg, R.~W., Hounsgaard, J. and Ditlevsen, S.} (2011).
\newblock Motoneuron membrane potentials follow a time inhomogeneous jump
  diffusion process.
\newblock {\em J. Comput. Neurosci.\/} To appear, available online.

\bibitem{Kurtz1978}
{\sc Kurtz, T.~G.} (1978).
\newblock Strong approximation theorems for density dependent {M}arkov chains.
\newblock {\em Stoch. Proc. Appl.\/} {\bf 6,} 223--240.

\bibitem{KuskeGordilloGreenwood2007}
{\sc Kuske, R., Gordillo, L.~F. and Greenwood, P.} (2007).
\newblock Sustained oscillations via coherence resonance in {S}{I}{R}.
\newblock {\em J. Theor. Biol.\/} {\bf 245,} 459--469.

\bibitem{LanskyDitlevsen2008}
{\sc Lansky, P. and Ditlevsen, S.} (2008).
\newblock A review of the methods for signal estimation in stochastic diffusion
  leaky integrate-and-fire neuronal models.
\newblock {\em Biol. Cybern.\/} {\bf 99,} 253--262.

\bibitem{Lapicque1907}
{\sc Lapicque, L.} (1907).
\newblock Recherches quantitatives sur l'excitation electrique des nerfs
  traitee comme une polarization.
\newblock {\em J. Physiol. Pathol. Gen.\/} {\bf 9,} 620--35.

\bibitem{MorrisLecar1981}
{\sc Morris, C. and Lecar, H.} (1981).
\newblock Voltage oscillations in the barnacle giant muscle fiber.
\newblock {\em Biophysics Journal\/} {\bf 35,} 193--213.

\bibitem{ISI:000237628500003}
{\sc Pfister, J.~P., Toyoizumi, T., Barber, D. and Gerstner, W.} ({2006}).
\newblock {Optimal spike-timing-dependent plasticity for precise action
  potential firing in supervised learning}.
\newblock {\em {Neural Comput.}\/} {\bf {18},} {1318--1348}.

\bibitem{RinzelErmentrout1998}
{\sc Rinzel, J. and Ermentrout, B.} (1998).
\newblock {\em In: Methods in neuronal modeling} 2nd~ed.
\newblock MIT Press, Cambridge.
\newblock ch.~Analysis of neural excitability and oscillations, pp.~251--291.

\bibitem{TatenoPakdaman2004}
{\sc Tateno, T. and Pakdaman, K.} (2004).
\newblock Random dynamics of the {M}orris-{L}ecar neural model.
\newblock {\em Chaos\/} {\bf 14,} 511--530.

\end{thebibliography}

\end{document}